\documentclass[12pt]{amsart}
\usepackage{graphicx}

\theoremstyle{plain}
\newtheorem{theorem}{\indent\sc Theorem}[section]
\newtheorem{lemma}[theorem]{\indent\sc Lemma}

\newtheorem{proposition}[theorem]{\indent\sc Proposition}

\theoremstyle{definition}

\newtheorem*{remark0}{\indent\sc Remark}

\title[Metrical theory for $\alpha$-Rosen fractions]{Metrical theory
for $\alpha$-Rosen fractions}

\author{Karma Dajani}
\address{
Department of Mathematics\endgraf Utrecht University\endgraf
Postbus 80.000\endgraf 3508 TA Utrecht\endgraf the Netherlands}
\email{dajani@math.uu.nl}
\author{Cor Kraaikamp}
\address{
TU Delft\endgraf EWI (DIAM)\endgraf Mekelweg 4\endgraf 2628 CD
Delft\endgraf the Netherlands} \email{c.kraaikamp@tudelft.nl}
\author{Wolfgang Steiner}
\address{
LIAFA, CNRS UMR 7089\endgraf Universit\'e Paris Diderot -- Paris
7\endgraf Case 7014\endgraf 75205 Paris Cedex 13\endgraf France}
\email{steiner@liafa.jussieu.fr}

\thanks{The third author was supported by NWO Bezoekersbeurs
{\bf NB 61-572}, by the Austrian Science Fund FWF, grant S8302-MAT,
and by the Agence Nationale de la Recherche, grant 
ANR-06-JCJC-0073.}

\begin{document}
\begin{abstract}
The Rosen fractions form an infinite family which generalizes the
nearest-integer continued fractions. In this paper we introduce a
new class of continued fractions related to the Rosen fractions,
the $\alpha$-Rosen fractions. The metrical properties of these
$\alpha$-Rosen fractions are studied.

We find planar natural extensions for the associated interval
maps, and show that their domains of definition are closely
related to the domains of the `classical' Rosen fractions. This
unifies and generalizes results of diophantine approximation from
the literature.
\end{abstract}

\maketitle

\section{Introduction}
Although David Rosen \cite{[Ros]} introduced as early as 1954 an
infinite family of continued fractions which generalize the
nearest-integer continued fraction, it is only very recently that
the metrical properties of these so-called Rosen fractions haven
been investigated; see e.g.~\cite{[Schm]}, \cite{[N2]}, \cite{[GH]}
and~\cite{[BKS]}. In this paper we will introduce $\alpha$-Rosen
fractions, and study their metrical properties for special choices
of~$\alpha$. These choices resemble Nakada's $\alpha$-expansions,
in fact for $q=3$ these are Nakada's $\alpha$-expansions; see
also~\cite{[N1]}. To be more precise, let $q\in\mathbb Z$, $q\geq
3$, and $\lambda = \lambda_q=2 \cos \frac{\pi}{q}$. Then we define
for $\alpha \in [\frac{1}{2}, \frac{1}{\lambda}]$ the map
$T_{\alpha} : [\lambda (\alpha -1), \lambda\alpha] \to [\lambda
(\alpha -1), \lambda\alpha)$ by
\begin{equation}\label{alfamap}
T_{\alpha}(x) := \left| \frac{1}{x}\right| -\lambda \left\lfloor
\left| \frac{1}{x\lambda}\right| + 1-\alpha \right\rfloor ,\,
x\neq 0,
\end{equation}
and $T_{\alpha}(0) :=0$. Here, $\lfloor \xi \rfloor$ denotes the
{\em floor} (or \emph{entier}) of $\xi$, i.e., the greatest
integer smaller than or equal to $\xi$. In order to have positive
digits, we demand that $\alpha \leq 1/\lambda$. Setting
$d(x)=\left\lfloor\left|\frac1{x\lambda}\right|+1-\alpha\right\rfloor$
(with $d(0)=\infty$), $\varepsilon(x)=\mathrm{sgn}(x)$, and more
generally
\begin{equation}\label{epsd}
\varepsilon_n(x)=\varepsilon_n =\varepsilon \left(
T_{\alpha}^{n-1}(x)\right) \quad \text{and}\quad
d_n(x)=d_n=d\left( T_{\alpha}^{n-1}(x)\right)
\end{equation}
for $n\geq 1$, one obtains for $x\in I_{q,\alpha}:=[\lambda (\alpha
-1),\alpha \lambda ]$ an expression of the form
$$
x = \frac{\varepsilon_1}{d_1 \lambda + {\displaystyle
\frac{\varepsilon_2}{d_2 \lambda + \cdots
+\cfrac{\varepsilon_n}{d_n\lambda + T^n_{\alpha}(x)}}}} ,
$$
where $\varepsilon_i \in \{\pm 1,0\}$ and $d_i \in \mathbb
N\cup\{\infty\}$. Setting
\begin{equation}\label{convergents}
\frac{R_n}{S_n}=\frac{\varepsilon_1}{d_1 \lambda + {\displaystyle
\frac{\varepsilon_2}{d_2 \lambda + \cdots
+\cfrac{\varepsilon_n}{d_n\lambda}}}} =: [\, \varepsilon_1:d_1, \,
\varepsilon_2:d_2, \dots , \varepsilon_n:d_n\,],
\end{equation}
we will show in Section~\ref{conv} that
$$
\lim_{n\to\infty} \frac{R_n}{S_n}=x,
$$
and for convenience we will write
\begin{equation}\label{cf}
x = \frac{\varepsilon_1}{d_1 \lambda + {\displaystyle
\frac{\varepsilon_2}{d_2 \lambda + \cdots}}} =: [\,
\varepsilon_1:d_1, \, \varepsilon_2:d_2, \dots\,].
\end{equation}
We call $R_n/S_n$ the $n$th $\alpha$-Rosen convergent of $x$,
and~(\ref{cf}) the $\alpha$-Rosen fraction of $x$.\smallskip\

The case $\alpha = 1/2$ yields the Rosen fractions, while the case
$\alpha = 1/\lambda$ is the Rosen fraction equivalent of the
classical regular continued fraction expansion (RCF). In case
$q=3$ (and $1/2\leq \alpha \leq 1/\lambda$), the above defined
$\alpha$-Rosen fractions are in fact Nakada's $\alpha$-expansions
(and the case $\alpha=1/\lambda =1$ is the RCF). Already from
\cite{[BKS]} it is clear that in order to construct the underlying
ergodic system for any $\alpha$-Rosen fraction and the planar
natural extension for the associated interval map $T_{\alpha}$, it
is fundamental to understand the orbit under $T_{\alpha}$ of the
two endpoints $\lambda (\alpha -1)$ and $\lambda\alpha$ of
$X=X_{\alpha} :=[\lambda (\alpha -1), \lambda\alpha]$. Although
the situation is in general more complicated than the `classical
case' from~\cite{[BKS]}, the natural extension together with the
invariant measure can be given, and it is shown that this
dynamical system is weakly Bernoulli.

Using the natural extension, metrical properties of the
$\alpha$-Rosen fractions will be given in Section~\ref{Metrical
properties}.

\section{Natural extensions}\label{General alfa's}
In this section we find the ``smallest'' domain
$\Omega_{\alpha}\subset \mathbb R^2$ on which the map
\begin{equation}\label{natural extension map}
\mathcal{T}_{\alpha} (x,y)=\left(
T_{\alpha}(x),\frac{1}{d(x)\lambda + \varepsilon(x)y}\right)
,\quad (x,y)\in \Omega_{\alpha},
\end{equation}
is bijective a.e.. We will deal with the general case, resembling
Nakada's $\alpha$-expansions, i.e., $1/2\le\alpha\le 1/\lambda$ and
$\lambda=\lambda_q=2\cos\pi/q$ for some fixed $q\in\mathbb Z$,
$q\geq 4$ (the case $q=3$ is in fact the case of Nakada's
$\alpha$-expansions; see also~\cite{[N1]}). As in~\cite{[BKS]}, we
need to discern between odd and even $q$'s, but some properties are
shared by both cases, and these are collected here first.\smallskip

For $x\in [\lambda (\alpha -1),\lambda \alpha]$, setting
$$
A_i=\left( \begin{array}{cc} 0 & \varepsilon_i\\
1 &  d_i\lambda
\end{array} \right) ,\quad \text{and}\,\,\,\, M_n=A_1\cdots A_n=
\left( \begin{array}{cc} K_n  & R_n\\
L_n & R_n
\end{array} \right) ,
$$
it immediately follows from $M_n=M_{n-1}A_n$ that $K_n=R_{n-1}$,
$L_n=S_{n-1}$, and
\begin{equation}\label{recurrence relations}
\begin{array}{cccl}
R_{-1}:=1, & R_0:=0, & R_n=d_n\lambda
R_{n-1}+\varepsilon_nR_{n-2}, &
\text{for }n=1,2,\dots\,, \vspace{1mm} \\
S_{-1}:=0, & S_0:=1, & S_n=d_n\lambda
S_{n-1}+\varepsilon_nS_{n-2}, & \text{for }n=1,2,\dots\,,
\end{array}
\end{equation}
if $d_n<\infty$. For a matrix $A=\left( \begin{array}{cc} a & b\\
c & d \end{array}\right)$, with $\det (A)\neq 0$, we define the
corresponding M\"obius (or fractional linear) transformation by
$$
A(x)=\frac{ax+b}{cx+d}.
$$
Consequently, considering $M_n$ as a M\"obius transformation, we
find that
$$
M_n(0)=\frac{R_n}{S_n},\quad \text{and}\,\, M_n(0) = A_1\cdots
A_n(0)=\dots = [\, \varepsilon_1:d_1, \, \varepsilon_2:d_2, \dots ,
\varepsilon_n:d_n\,].
$$
It follows that the numerators and denominators of the
$\alpha$-Rosen convergents of $x$ from~(\ref{convergents}) satisfy
the usual recurrence relations~(\ref{recurrence relations}); see
also \cite{[BKS]}, p.~1279.

Furthermore, since
$$
x=M_{n-1}\left( \begin{array}{cc} 0 & \varepsilon_n\\
1 & d_n\lambda +T_{\alpha}^n(x)
\end{array}\right) (0),
$$
we have that
\begin{equation}\label{xvsTn}
x=\frac{R_n+T_\alpha^n(x)R_{n-1}}{S_n+T_\alpha^n(x)S_{n-1}}\
\text{ and }\ T_\alpha^n(x)=\frac{R_n-S_nx}{S_{n-1}x-R_{n-1}}.
\end{equation}

Let $\ell_0=(\alpha-1)\lambda$ be the left-endpoint of the interval
on which the continued fraction map $T_\alpha$ was defined in
(\ref{alfamap}), $r_0=\alpha\lambda$ its right-endpoint and let
$$
\Delta(\varepsilon:d) =
\{\,x\in[(\alpha-1)\lambda,\alpha\lambda]\mid
\varepsilon_1(x)=\varepsilon,\ d_1(x)=d\,\},
$$
be the cylinders of order $1$ of numbers with same first digits
given by (\ref{cf}). If we set
$$
\delta_d=\frac1{(\alpha+d)\lambda}
$$
for all $d\ge 1$, then the cylinders are given by the following
table
$$
\begin{array}{c|c|c|c|c}
\Delta(-1:1) & \Delta(-1:d),\,d\ge2 & \Delta(0:\infty) &
\Delta(+1:d),\,d\ge2 & \Delta(+1:1) \\
\hline [\ell_0,\delta_1) & [-\delta_{d-1},-\delta_d) & \{0\} &
(\delta_d,\delta_{d-1}] & (\delta_1,r_0]
\end{array}
$$
where we have used that $r_0>\delta_1$ since $\lambda\ge\sqrt 2$ for
$q\geq 4$. Note that we have by definition that
$$
T_\alpha(x)=\varepsilon/x-\lambda d,
$$
for all $x\in\Delta(\varepsilon,d)$, $x\ne 0$.

\begin{figure}[h!t]
\includegraphics[scale=.99]{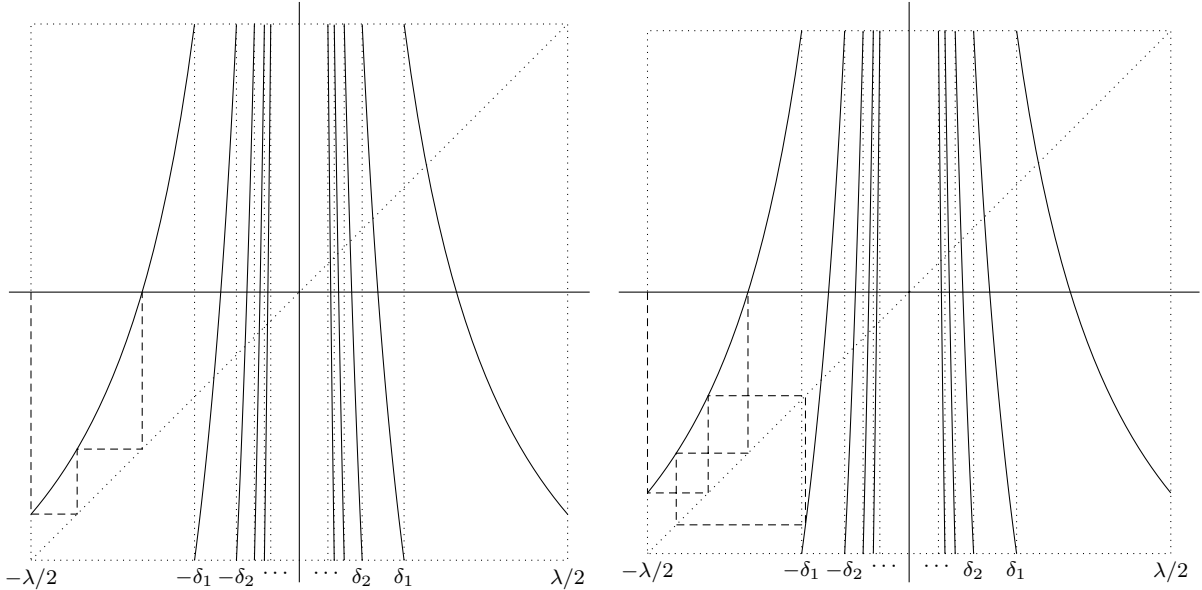}
\caption{The map $T_{1/2}$ and the orbit of $-\lambda/2$ (dashed broken 
line) for $q=8$ (left) and $q=7$ (right).} \label{fig:phi} 
\end{figure}

Setting $\ell_n=T^n_\alpha(\ell_0)$, $r_n=T^n_\alpha(r_0)$, $n\geq
0$, we have that
$$
r_1=\frac{1}{\alpha\lambda}-\lambda=
-\frac{\alpha\lambda^2-1}{\alpha\lambda}<0.
$$
In case $\alpha=1/2$, we write $\phi_n$ instead of $\ell_n$, for
$n\geq 0$. In~\cite{[BKS]}, it was shown that
$$
-\lambda /2= \begin{cases}
[ \, (-1\, :\, 1)^{p-1} \, ] , &  \text{if $q=2p$}\\
\text{[}\, (-1\, :\, 1)^h, \, -1\, :\, 2, \, ( -1 \, :\, 1 )^h \,
\text{]} \, , & \text{if $q=2h+3$},
\end{cases}
$$
from which it immediately follows that
\begin{equation}\label{evenphi's}
\phi_0=-\frac\lambda2<\phi_1<\cdots<\phi_{p-2}=-\frac1\lambda<\phi_{p-1}
=0, \quad \text{if $q=2p$},
\end{equation}
and that for $q=2h+3$,
$$
\phi_0=-\frac\lambda2<\phi_1<\cdots<\phi_{h-2}<\phi_{h-1}<-\frac2{3\lambda}<
\phi_h<-\frac2{5\lambda},
$$
$\phi_0<\phi_{h+1}<\phi_1$, $\phi_{h+1}=1-\lambda$, and
$$
\phi_{h+1}<\phi_{h+2}<\cdots<\phi_{2h}=-\frac1\lambda<-\frac2{3\lambda}
<\phi_{2h+1}=0;
$$
see also Figure~\ref{fig:phi}.

Thus we see that the behavior of the orbit of $-\lambda /2$ is
very different in the even case compared to the odd case; see also
Figure~\ref{fig:situationschets}, where the relevant terms of
$(\phi_n)_{n\geq 0}$, $(\ell_n)_{n\geq 0}$, and $(r_n)_{n\geq 0}$
are displayed for even $q$.\medskip\

Direct verification yields the following lemma.
\begin{lemma}\label{lemma1}
For $q\geq 4$, $1/2\le\alpha\le1/\lambda$, we have that
$$
\phi_0=-\frac\lambda2\le\ell_0\le
r_1<\phi_1=-\frac{\lambda^2-2}\lambda,
$$
with $\phi_0=\ell_0$ if and only if $\alpha=1/2$, and $\ell_0=r_1$
if and only if $\alpha=1/\lambda$.
\end{lemma}

In~\cite{[BKS]}, the sequence $(\phi_n)_{n\geq 0}$ plays a crucial
role in the construction of the natural extension of the Rosen
fractions. Due to the fact that the orbits of both $-\lambda/2$
and $\lambda/2$ would become constant $0$ after a finite number of
steps (depending on $q$), the natural extension of the Rosen
fraction could be easily constructed. In this paper, the
$(\ell_n)_{n\geq 0}$ and $(r_n)_{n\geq 0}$ play a role comparable
to that of the sequence $(\phi_n)_{n\geq 0}$ (even though the
$\phi_n$'s are frequently used as well).

Let $x\in [\lambda (\alpha -1),\alpha \lambda ]$ be such, that
$(\varepsilon_n(x):d_n(x))=(-1:1)$ for $n=1,2,\ldots,m$. Then it
follows from~(\ref{recurrence relations}) that the $\alpha$-Rosen
convergents of $x$ satisfy
\begin{eqnarray*}
R_{-1}=1, & R_0=0, & R_n=\lambda R_{n-1}-R_{n-2},\quad
\text{for }n=1,2,\dots,m\,,\\
S_{-1}=0, & S_0=1, & S_n=\lambda S_{n-1}-S_{n-2},\quad \text{for
}n=1,2,\dots,m\,.
\end{eqnarray*}
As in~\cite{[BKS]}, we define the auxiliary sequence
$(B_n)_{n\geq0}$ by
\begin{equation}\label{Bn's}
B_0 = 0,\ B_1 = 1,\quad B_n = \lambda B_{n-1} - B_{n-2},\quad
\text{for } n=2,3,\ldots \, .
\end{equation}
This yields for $n=1,\dots,m$ that $R_n=-B_n$, $S_n=B_{n+1}$, and
$T_\alpha^n(x)=-\frac{B_n+B_{n+1}x}{B_{n-1}+B_nx}$ by
(\ref{xvsTn}). It follows that
\begin{equation}\label{ln}
\ell_n =
-\frac{B_n+B_{n+1}(\alpha-1)\lambda}{B_{n-1}+B_n(\alpha-1)\lambda}
=
-\frac{B_{n+1}\alpha\lambda-B_{n+2}}{B_n\alpha\lambda-B_{n+1}}\quad
\mbox{ if  }\ell_0=[(-1:1)^n,\ldots\,].
\end{equation}

For $x=[+1:1,(-1:1)^{n-1},\ldots\,]=-[(-1:1)^n,\ldots\,]\,$, we
obtain similarly $R_n=B_n$, $S_n=B_{n+1}$, thus
$T_\alpha^n(x)=\frac{B_n-B_{n+1}x}{B_nx-B_{n-1}}$ and
\begin{equation}\label{rn}
r_n = -\frac{B_{n+1}\alpha\lambda-B_n}{B_n\alpha\lambda-B_{n-1}}
\qquad \text{if } r_0=[+1:1,(-1:1)^{n-1},\ldots\,]\,.
\end{equation}
It is easy to see that $B_n=\sin\frac{n\pi}q/\sin\frac\pi q$, hence
$(B_n)_{n\ge0}$ is a periodic sequence with period length $2q$.

\subsection{Even indices} Let $q=2p$, $p\in\mathbb N$, $p\geq 2$.
Essential in the construction of the natural extension is the
following theorem.
\begin{theorem}\label{thm:evenorder}
Let $q=2p$, $p\in\mathbb N$, $p\geq 2$, and let the sequences
$(\ell_n)_{n\geq 0}$ and $(r_n)_{n\geq 0}$ be defined as before.
If $1/2<\alpha<1/\lambda$, then we have that
\begin{equation}\label{evenorder}
\ell_0<r_1<\ell_1<\cdots<r_{p-2}<\ell_{p-2}<-\delta_1<r_{p-1}<0<
\ell_{p-1}<r_0,
\end{equation}
$d_p(r_0)=d_p(\ell_0)+1$ and $\ell_p=r_p$. If $\alpha=1/2$, then
we have that
\begin{equation}\label{evenorder 1/2}
\ell_0<r_1=\ell_1<\cdots<r_{p-2}=\ell_{p-2}<-\delta_1<r_{p-1}=
\ell_{p-1}=0<r_0.
\end{equation}
If $\alpha=1/\lambda$, then we have that
\begin{equation}\label{evenorder 1/lambda}
\ell_0=r_1<\ell_1=r_2<\cdots<\ell_{p-2}=r_{p-1}=-\delta_1<0<r_0.
\end{equation}
\end{theorem}

\begin{proof}
If $\alpha = 1/2$, then $\ell_0=\phi_0$ and $r_0=-\phi_0$,
hence~(\ref{evenorder 1/2}) is an immediate consequence
of~(\ref{evenphi's}).

In general, in view of Lemma~\ref{lemma1} and the fact that
$\phi_0=[(-1:1)^{p-1}]$, we have the following situation:
$T_\alpha([\ell_0,\phi_1))=[\ell_1,\phi_2)$ and
$T_\alpha([\phi_{j-1},\phi_j))=[\phi_j,\phi_{j+1})$ for
$j=2,3,\ldots,p-2$, cf.\ Figure~\ref{fig:situationschets}. This
yields that $\ell_0=[(-1:1)^{p-2},\ldots\,].$

\begin{figure}[h!t]
\includegraphics{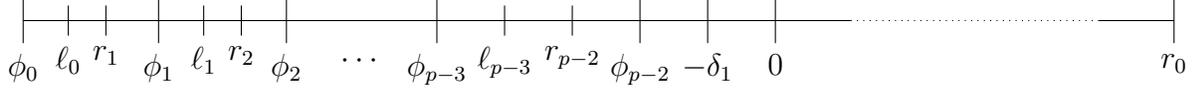}
\caption[situation]{\label{fig:situationschets} The relevant terms
of $(\phi_n)_{n\geq 0}$, $(\ell_n)_{n\geq 0}$, and $(r_n)_{n\geq
0}$ for even $q$.}
\end{figure}

Since $\sin\frac{(p-1)\pi}{2p}=\sin\frac{(p+1)\pi}{2p}$, we obtain
$$
B_{p-1}=B_{p+1},\quad B_{p-1}=\frac\lambda2 B_p,\quad
B_{p-2}=\left(\frac{\lambda^2}2-1\right)B_p.
$$
By (\ref{ln}), we have therefore that
$$
\ell_{p-2} =
-\frac{B_{p-1}\alpha\lambda-B_p}{B_{p-2}\alpha\lambda-B_{p-1}} =
-\frac{2-\alpha\lambda^2}{\lambda(1-\alpha\lambda^2+2\alpha)} \le
-\frac1{(\alpha+1)\lambda} = -\delta_1,
$$
with $\ell_{p-2}=-\delta_1$ if and only if $\alpha=1/\lambda$. For
$\alpha=1/\lambda$, we clearly have that $r_0=1$, thus
$r_1=1-\lambda=\ell_0$ and (\ref{evenorder 1/lambda}) is proved.

If $1/2<\alpha<1/\lambda$, then we have that $\ell_{p-2}<-\delta_1$,
hence $\ell_0=[(-1:1)^{p-1},\ldots , (1:d_p),\dots \,]$, with
$d_p\geq 1$, and again due to (\ref{ln}) we obtain
$$
\ell_{p-1} =
-\frac{B_p\alpha\lambda-B_{p-1}}{B_{p-1}\alpha\lambda-B_p} =
\frac{(2\alpha-1)\lambda}{2-\alpha\lambda^2}>0.
$$
Similarly, we have that $r_0=[+1:1,(-1:1)^{p-2},\ldots\,]$ and, by
(\ref{rn}),
$$
r_{p-1} =
-\frac{B_p\alpha\lambda-B_{p-1}}{B_{p-1}\alpha\lambda-B_{p-2}} =
-\frac{(2\alpha-1)\lambda}{2-(1-\alpha)\lambda^2}
\in(-\delta_1,0),
$$
thus (\ref{evenorder}) is proved. Since
$$
\left|\frac1{r_{p-1}}\right|-\left|\frac1{\ell_{p-1}}\right|=\lambda,
$$
it follows from the definition of $T_\alpha$ that $\ell_p=r_p$ and
$d_1(r_{p-1})=d_1(\ell_{p-1})+1$. With
$d_p(x)=d_1(T_\alpha^{p-1}(x))$, the theorem is proved.
\end{proof}

\begin{remark0}
The structure of the $\ell_n$'s and $r_n$'s allows us to determine
all possible sequences of ``digits''. For example, the longest
consecutive sequence of digits $(-1:1)$ contains $p-1$ terms if
$\alpha<1/\lambda$ since $\ell_{p-2}<-\delta_1$ and
$\ell_{p-1}\ge-\delta_1$. In case $\alpha=1/\lambda$, we only have
$(-1:1)^{p-2}$. In particular in case $q=4$, $\alpha=1/\lambda$, the
cylinder $\Delta(-1:1)$ is empty.

On the other hand, $(+1:1)$ is always followed by $(-1:1)^{p-2}$
since $r_{p-2}<-\delta_1$, with $\Delta(+1:1)=\{r_0\}$ in case
$\alpha=1/\lambda$.
\end{remark0}

Now we construct the domain $\Omega_\alpha$ upon which $\mathcal
T_\alpha$ is bijective.

\begin{theorem}\label{thm:even case nat ext}
Let $q=2p$ with $p\ge 2$. Then the system of relations
$$
\left\{\begin{array}{rl}
(\mathcal R_1): & \qquad H_1=1/(\lambda+H_{2p-1}) \\
(\mathcal R_2): & \qquad H_2=1/\lambda \\
(\mathcal R_n): & \qquad H_n=1/(\lambda-H_{n-2})\qquad
\text{for }n=3,4,\dots,2p-1 \\
(\mathcal R_{2p}): & \qquad H_{2p-2}=\lambda/2 \\
(\mathcal R_{2p+1}): & \qquad H_{2p-3}+H_{2p-1}=\lambda \\
\end{array}\right.
$$
admits the (unique) solution
\begin{gather*}
H_{2n} = -\phi_{p-n} = \frac{B_n}{B_{n+1}} =
\frac{\sin\frac{n\pi}{2p}}{\sin\frac{(n+1)\pi}{2p}}\ \text{ for }
n=1,2,\ldots,p-1, \\
H_{2n-1} = \frac{B_{p-n}-B_{p+1-n}}{B_{p-1-n}-B_{p-n}} =
\frac{\cos\frac{n\pi}{2p}-\cos\frac{(n-1)\pi}{2p}}
{\cos\frac{(n+1)\pi}{2p}-\cos\frac{n\pi}{2p}}\ \text{ for }
n=1,2,\ldots,p,
\end{gather*}
in particular $H_{2p-1}=1$.

Let $1/2<\alpha<1/\lambda$ and
$\Omega_\alpha=\bigcup_{n=1}^{2p-1}J_n\times[0,H_n]$ with
$J_{2n-1}=[\ell_{n-1},r_n)$, $J_{2n}=[r_n,\ell_n)$ for
$n=1,2,\ldots,p-1$, and $J_{2p-1}=[\ell_{p-1},r_0)$. Then the map
$\mathcal T_\alpha:\, \Omega_\alpha \rightarrow \Omega_\alpha$
given by (\ref{natural extension map}) is bijective off of a set
of Lebesgue measure zero.
\end{theorem}

\begin{figure}[h!t]
\includegraphics{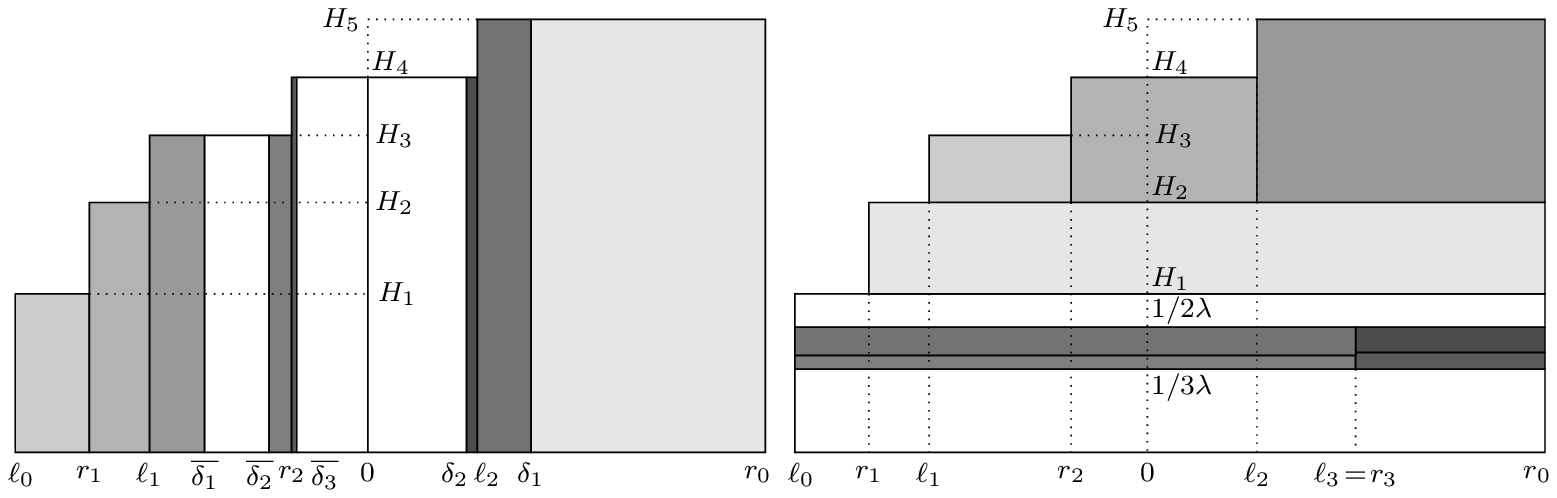}
\caption[natural extension q=6]{\label{fig:nat ext q=6} The
natural extension domain $\Omega_\alpha$ (left) and its image
under $\mathcal T_\alpha$ (right) of the $\alpha$-Rosen continued
fraction ($\overline{\delta_n}=-\delta_n$); here $q=6$,
$\alpha=0.53$, $d_p(\ell_0)=2$, $d_p(r_0)=3$.}
\end{figure}

\begin{proof}
It is  easily seen that the solution of this system of relations
is unique and valid, and that $\mathcal T_\alpha$ is injective. We
thus concern ourselves with the surjectivity of $\mathcal
T_\alpha$; see also Figure~\ref{fig:nat ext q=6}.

By (\ref{evenorder}), we have $J_{n-2}\subset \Delta(-1:1)$ for
$n=3,4,\ldots,2p-2$, thus
$$
\mathcal T_\alpha(J_{n-2}\times[0,H_{n-2}]) =
J_n\times\left[\frac1\lambda,\frac1{\lambda-H_{n-2}}\right] =
J_n\times[H_2,H_n]\,,
$$
where we have used ($\mathcal R_2$) and ($\mathcal R_n$).
Furthermore, $(\mathcal R_{2p-1})$ gives
$$
\mathcal T_\alpha([\ell_{p-2},-\delta_1)\times[0,H_{2p-3}]) =
[\ell_{p-1},r_0)\times\left[\frac1\lambda,\frac1{\lambda-H_{2p-3}}
\right] = J_{2p-1}\times[H_2,H_{2p-1}]\,.
$$
For $n=2,3,\ldots,d_1(r_{p-1})-1=d_p(r_0)-1$, we have that
$$
\mathcal
T_\alpha([-\delta_{n-1},-\delta_n)\times[0,H_{2p-3}])=[\ell_0,
r_0)\times\left[\frac1{n\lambda},\frac1{n\lambda-H_{2p-3}}\right].
$$
The remaining part of the rectangle $J_{2p-3}\times[0,H_{2p-3}]$
is mapped to
$$
\mathcal
T_\alpha([-\delta_{d_p(r_0)-1},r_{p-1})\times[0,H_{2p-3}])=
[\ell_0,r_p) \times
\left[\frac1{d_p(r_0)\lambda},\frac1{d_p(r_0)\lambda-H_{2p-3}}\right].
$$

Now consider the image of $J_{2p-1}\times[0,H_{2p-1}]$. If
$d_p(\ell_0)\geq2$, then it is split into
\begin{gather*}
\mathcal
T_\alpha((\ell_{p-1},\delta_{d_p(\ell_0)-1}]\times[0,H_{2p-1}] =
[\ell_0,r_p) \times \left[\frac1{d_p(\ell_0)\lambda+H_{2p-1}},
\frac1{d_p(\ell_0)\lambda}\right], \\
\mathcal T_\alpha((-\delta_n,-\delta_{n-1}]\times[0,H_{2p-1}]) =
[\ell_0,r_0)\times\left[\frac1{n\lambda+H_{2p-1}},\frac1{n\lambda}
\right]\ \text{ for }n=2,3,\ldots,d_p(\ell_0)-1, \\
\mathcal T_\alpha((\delta_1,r_0)\times[0,H_{2p-1}]) =
(r_1,r_0)\times
\left[\frac1{\lambda+H_{2p-1}},\frac1\lambda\right] =
(r_1,r_0)\times[H_1,H_2],
\end{gather*}
where we have used ($\mathcal R_1$). Since
$H_{2p-3}+H_{2p-1}=\lambda$ and $d_p(r_0)=d_p(\ell_0)+1$, the
different parts of $\mathcal
T_\alpha([-\delta_1,r_{p-1})\times[0,H_{2p-3}])$ and $\mathcal
T_\alpha((\ell_{p-1},\delta_1]\times[0,H_{2p-1}])$ ``layer one
under the other'' and ``fill up like a jig-saw puzzle''
$$
\left([\ell_0,r_p)\times
\left[\frac1{d_p(r_0)\lambda},\frac1{d_p(\ell_0)\lambda}\right]\right)
\cup
\left([\ell_0,r_0)\times\left[\frac1{d_p(r_0)\lambda},H_1\right]\right).
$$
In case $d_p(\ell_0)=1$, we simply have
\begin{gather*}
\mathcal T_\alpha([-\delta_1,r_{p-1})\times[0,H_{2p-3}])=
[\ell_0,r_p)\times[1/(2\lambda),H_1]\,, \\
\mathcal T_\alpha((\ell_{p-1},r_0)\times[0,H_{2p-1}] =
(r_1,r_p)\times[H_1,H_2]\,.
\end{gather*}

Finally, the image of the central rectangle
$J_{2p-2}\times[0,H_{2p-2}]$ is split into
\begin{gather*}
\mathcal T_\alpha([r_{p-1},-\delta_{d_p(r_0)})\times[0,H_{2p-2}])
= [r_p,r_0) \times
\left[\frac1{d_p(r_0)\lambda},\frac1{d_p(r_0)\lambda-H_{2p-2}}\right],
\\ \mathcal T_\alpha([-\delta_{n-1},-\delta_n)\times[0,H_{2p-2}]) =
[\ell_0,r_0)\times\left[\frac1{n\lambda},\frac1{n\lambda-H_{2p-2}}
\right]\ \text{ for }n>d_p(r_0), \\
\mathcal T_\alpha((\delta_n,\delta_{n-1}]\times[0,H_{2p-2}]) =
[\ell_0,r_0)\times\left[\frac1{n\lambda+H_{2p-2}},\frac1{n\lambda}
\right]\ \text{ for }n>d_p(\ell_0), \\
\mathcal
T_\alpha((\delta_{d_p(\ell_0)},\ell_{p-1}]\times[0,H_{2p-2}]) =
[r_p,r_0) \times \left[\frac1{d_p(\ell_0)\lambda+H_{2p-2}},
\frac1{d_p(\ell_0)\lambda}\right].
\end{gather*}
Since $H_{2p-2}=\lambda/2$ and $d_p(r_0)=d_p(\ell_0)+1$, the union
of these images is
$$
\left([\ell_0,r_0)\times\left(0,\frac1{d_p(r_0)\lambda}\right]\right)
\cup \left([r_p,r_0)\times
\left[\frac1{d_p(r_0)\lambda},\frac1{d_p(\ell_0)\lambda}\right]\right).
$$
Therefore $\mathcal T_\alpha(\Omega_\alpha)$ and $\Omega_\alpha$
differ only by a set of Lebesgue measure zero.
\end{proof}

\begin{remark0}
If $\alpha=1/2$, then the intervals $J_{2n}$ are empty and
$r_{p-1}=\ell_{p-1}=0$. The proof of Theorem~\ref{thm:even case
nat ext} remains valid, with
$d_1(r_{p-1})=d_1(\ell_{p-1})=\infty$; see also~\cite{[BKS]}.
Since $\ell_n=\phi_n$ for $n=0,1,\ldots,p-1$, we have
$$
\Omega_{1/2} = \bigcup_{n=1}^{p-1}
\left(\left[-\frac{B_n}{B_{n+1}},-\frac{B_{n-1}}{B_n}\right)
\times \left[0,\frac{B_n-B_{n+1}}{B_{n-1}-B_n}\right]\right) \cup
\left(\left[0,\frac\lambda2\right)\times[0,1]\right).
$$
\end{remark0}

For $\alpha=1/\lambda$, we just have the intervals $J_{2n}$,
$n=1,2,\ldots,p-2$ and add
$J_{2p-2}=[r_{p-1},r_0)(=[-\delta_1,1))$. Furthermore, we have
that $r_n=\frac{B_n-B_{n+1}}{B_n-B_{n-1}}$ for $n=1,2,\ldots,p-1$
and $r_0=\frac{B_p-B_{p+1}}{B_p-B_{p-1}}=1$. This provides the
following theorem.
\begin{theorem}\label{thm:even case nat ext 1/lambda}
Let $q=2p$ with $p\ge 2$ and
$$
\Omega_{1/\lambda} = \bigcup_{n=1}^{p-1}
\left[\frac{B_n-B_{n+1}}{B_n-B_{n-1}},
\frac{B_{n+1}-B_{n+2}}{B_{n+1}-B_n}\right) \times
\left[0,\frac{B_n}{B_{n+1}}\right].
$$
Then $\mathcal T_{1/\lambda}:\, \Omega_{1/\lambda} \rightarrow
\Omega_{1/\lambda}$ is bijective off of a set of Lebesgue measure
zero.
\end{theorem}

\begin{proof}
By (\ref{evenorder 1/lambda}), we have that
$\ell_0=r_1<\ell_1=r_2<\ldots<\ell_{p-2}=r_{p-1}=-\delta_1$, thus
\begin{multline*}
\mathcal T_{1/\lambda}
\left(\left[\frac{B_n-B_{n+1}}{B_n-B_{n-1}},
\frac{B_{n+1}-B_{n+2}}{B_{n+1}-B_n}\right)\times
\left[0,\frac{B_n}{B_{n+1}}\right]\right) \\
= \left[\frac{B_{n+1}-B_{n+2}}{B_{n+1}-B_n},
\frac{B_{n+2}-B_{n+3}}{B_{n+2}-B_{n+1}}\right)\times
\left[\frac1\lambda,\frac{B_{n+1}}{B_{n+2}}\right] \quad\text{for
}n=1,2,\ldots,p-2.
\end{multline*}
The different parts of $[-r_{p-1},r_0)$ are mapped to
\begin{align*}
\mathcal T_{1/\lambda}
\left([-\delta_{n-1},-\delta_n)\times\left[0,\frac\lambda2\right]\right)
& =
[1-\lambda,1)\times\left[\frac1{n\lambda},\frac2{(2n-1)\lambda}\right]\
\text { for }n=2,3,\ldots \\
\mathcal T_{1/\lambda}
\left((\delta_n,\delta_{n-1}]\times\left[0,\frac\lambda2\right]\right)
& =
[1-\lambda,1)\times\left[\frac2{(2n+1)\lambda},\frac1{n\lambda}\right]\
\text{ for }n=2,3,\ldots \\
\mathcal T_{1/\lambda}
\left((\delta_1,1)\times\left[0,\frac\lambda2\right]\right) & =
[1-\lambda,1)\times\left[\frac2{3\lambda},\frac1\lambda\right]
\end{align*}
and the union of these images is
$[1-\lambda,1)\times(0,1/\lambda]$.
\end{proof}

\begin{remark0}
Note that there is a simple relation between $\Omega_{1/2}$ and
$\Omega_{1/\lambda}$, which will be useful in Section~\ref{mixing
properties}; reflect $\Omega_{1/2}$ in the line $y=x$ in case
$x\geq 0$, and reflect $\Omega_{1/2}$ in the line $y=-x$ in case
$x\leq 0$, to find $\Omega_{1/\lambda}$; see also
Figure~\ref{fig:nat ext 1/2,1/lambda}.
\end{remark0}

\begin{figure}[h!t]
\includegraphics{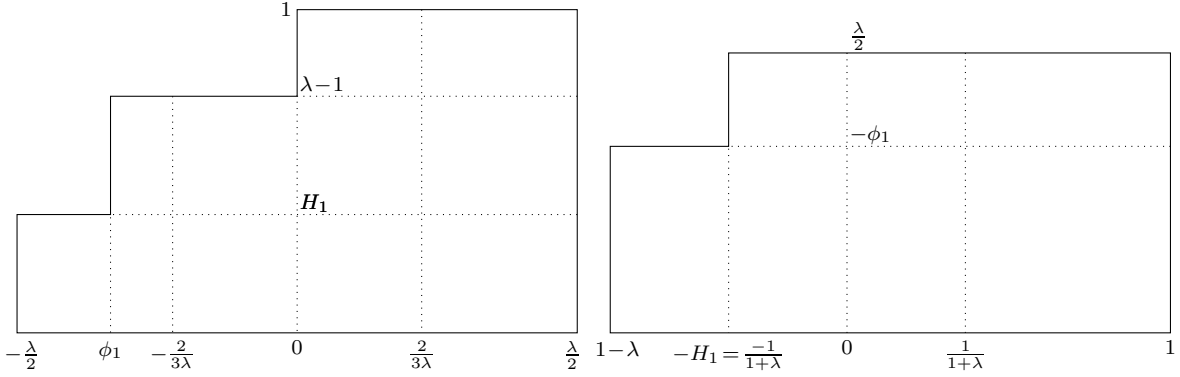}
\caption[1/2 and 1/lambda]{\label{fig:nat ext 1/2,1/lambda}
$\Omega_{1/2}$ (left) and $\Omega_{1/\lambda}$ (right); here
$q=6$.}
\end{figure}

As in~\cite{[BKS]}, a Jacobian calculation shows that $\mathcal
T_\alpha$ preserves the probability measure $\nu_\alpha$ with
density
$$
\frac{C_{q,\alpha}}{(1+xy)^2} \, ,
$$
where $C_{q,\alpha}$ is a normalizing constant. For the
calculation of this constant, we need the following lemma.
\begin{lemma} \label{lem:B-identity}
If $m_1-m_2=m_3-m_4$, then we have that
$$
B_{n+m_1}B_{-n+m_2}-B_{n+m_3}B_{-n+m_4}=B_{m_1-m_3}B_{m_2+m_3}
\quad \text{for all $n\in\mathbb Z$.}
$$
\end{lemma}

\begin{proof}
With $\zeta=\exp(\pi i/q)$, we have that
\begin{multline*}
B_{n+m_1}B_{-n+m_2}-B_{n+m_3}B_{-n+m_4} \\
=\frac{(\zeta^{n+m_1}-\zeta^{-n-m_1})(\zeta^{-n+m_2}-\zeta^{n-m_2})
-(\zeta^{n+m_3}-\zeta^{-n-m_3})(\zeta^{-n+m_4}-\zeta^{n-m_4})}
{(\zeta-\zeta^{-1})^2} \\ =
\frac{\zeta^{m_1+m_2}-\zeta^{-m_1-m_2}-\zeta^{m_3+m_4}-\zeta^{-m_3-m_4}}
{(\zeta-\zeta^{-1})^2} = B_{m_1-m_3}B_{m_2+m_3}.
\end{multline*}
\end{proof}

\begin{proposition}\label{lem:norm const even}
For $1/2\le\alpha\le1/\lambda$, the normalizing constant is
$$
C_{q,\alpha}=1/\log\frac{1+\cos\frac\pi q}{\sin\frac\pi q}.
$$
\end{proposition}

\begin{proof}
Similarly to~\cite{[BKS]}, integration of the density over
$\Omega_\alpha$ gives
$$
C_{q,\alpha}= 1/\log \left( \frac{1+r_0}{1+\ell_{p-1}}
\prod_{n=1}^{p-1} \frac{1+r_nH_{2n-1}}{1+\ell_{n-1}H_{2n-1}}
\frac{1+\ell_nH_{2n}}{1+r_nH_{2n}} \right)
$$
for $1/2<\alpha<1/\lambda$, by Theorem~\ref{thm:even case nat
ext}. Using (\ref{ln}), (\ref{rn}) and Lemma~\ref{lem:B-identity},
we find
$$
\frac{1+r_nH_{2n-1}}{1+\ell_{n-1}H_{2n-1}} =
\frac{B_n-B_{n-1}\alpha\lambda}{B_n\alpha\lambda-B_{n-1}},\qquad
\frac{1+\ell_nH_{2n}}{1+r_nH_{2n}} =
\frac{B_n\alpha\lambda-B_{n-1}}{B_{n+1}-B_n\alpha\lambda}
$$
for $n=1,2,\ldots,p-1$, and
$$
\frac{1+r_0}{1+\ell_{p-1}} =
\frac{B_p-B_{p-1}\alpha\lambda}{B_p-B_{p-1}}.
$$
Putting everything together, we obtain
$$
C_{q,\alpha} = 1/\log \frac1{B_p-B_{p-1}} = 1/\log
\frac{\sin\frac\pi q}{1-\cos\frac\pi q} = 1/\log
\frac{1+\cos\frac\pi q}{\sin\frac\pi q}.
$$
For $\alpha=1/2$, we have the same constant by the remark
following Theorem~\ref{thm:even case nat ext} and by \cite{[BKS]}.
Finally, the remark following Theorem~\ref{thm:even case nat ext
1/lambda} shows that $C_{q,1/\lambda}$ is the same constant as
well.
\end{proof}

Let $\mu_{\alpha}$ be the projection of $\nu_{\alpha}$ on the
first coordinate, let $\bar{\mathcal B}$ be the restriction of the
two-dimensional $\sigma$-algebra on $\Omega_{\alpha}$, and
${\mathcal B}$ be the Lebesgue $\sigma$-algebra on
$I_{q,\alpha}=[\lambda (\alpha -1),\alpha \lambda ]$.
In~\cite{[Roh]}, Rohlin introduced and studied the concept of
\emph{natural extension} of a dynamical system. In our setting, a
natural extension of $(I_{q,\alpha}, {\mathcal B},
\mu_{\alpha},T_{\alpha})$ is an invertible dynamical system
$(X_{\alpha}, {\mathcal B}_{X_{\alpha}}, \rho_{\alpha}, {\mathcal
S}_{\alpha})$, which contains $(I_{q,\alpha}, {\mathcal B},
\mu_{\alpha},T_{\alpha})$ as a factor, such that the Borel
$\sigma$-algebra ${\mathcal B}_{X_{\alpha}}$ of $X_{\alpha}$ is
the smallest ${\mathcal S}_{\alpha}$-invariant $\sigma$-algebra
that contains $\pi^{-1}({\mathcal B}_{X_{\alpha}})$, where $\pi$
is the factor map. A natural extension is unique up to
isomorphism.

We have the following theorem.
\begin{theorem}\label{thm:natural-extension-even-case}
Let $q\geq 4$, $q=2p$, and let $\frac{1}{2}\leq \alpha \leq
\frac{1}{\lambda}$. Then the dynamical system $(\Omega_{\alpha},
\bar{\mathcal B}, \nu_{\alpha},{\mathcal T}_{\alpha})$ is the
natural extension of the dynamical system $(I_{q,\alpha},
{\mathcal B}, \mu_{\alpha},T_{\alpha})$.
\end{theorem}

\begin{proof}
In case $\alpha =1/2$ and $\alpha =1/\lambda$ the proof is a
straightforward application of Theorem~21.2.2 from~\cite{[Schw]};
see also Examples~21.3.1 (the case of the RCF) and 21.3.2 (the
NICF) in~\cite{[Schw]}. However, for $1/2<\alpha <1/\lambda$, some
extra work is needed.\smallskip\

To be more precise, let us recall some terminology
from~\cite{[Schw]}. Let $B$ be a set, and $T: B\to B$ be a map.
The pair $(B,T)$ is called a \emph{fibred system} if the following
three conditions are satisfied:
\begin{itemize}
\item[(a)] There is a finite or countable set $I$ (called the
digit set);
\item[(b)] There is a map $k: B\to I$. Then the sets
$B(i)=k^{-1}(i)$ form a partition of $B$;
\item[(c)] The restriction of $B$ to any $B(i)$ is an injective
map;
\end{itemize}
see~\cite{[Schw]}, Definition~1.1.1. 
Clearly, in our setting 
$B=I_{q,\alpha}=[\lambda(\alpha-1),\alpha\lambda]$, $T=T_\alpha$, 
$I=\{\varepsilon:d;\;\varepsilon\in\{\pm 1,0\},
d\in\mathbb N\cup\{\infty\}\}$, and
$B(\varepsilon:d)=\Delta(\varepsilon:d)$.

The pair $(B^{\#},T^{\#})$ is called a \emph{dual fibred system}
(or \emph{backward algorithm}) with respect to $(B,T)$ if the
following condition holds: $(k_1,k_2,\dots,k_n)$ is an admissible
block of digits for $T$ if and only if $(k_n,\dots, k_2,k_1)$ is
admissible for $T^{\#}$; see Definition~21.1.1 in~\cite{[Schw]}.
Furthermore,
\begin{eqnarray*}
D(x)&:=&\{ y\in B^{\#};\, y\in B^{\#}(k_1,k_2,\dots,k_N) 
\text{ if and only if }\\
& & T^{-N}(x)\cap B(k_N,\dots,k_2,k_1)\neq \emptyset 
\text{ for all $N\geq 1$}\} ;
\end{eqnarray*}
see Definition~21.1.7 from~\cite{[Schw]}. The local inverse of the
map $T: B(k)\to B$ is denoted by $V(k)$. F.~Schweiger obtained the
following theorem; see~\cite{[Schw]}, Theorem~21.2.1.
\begin{theorem} Consider the following dynamical system
$(\bar{B},\bar{T})$, with $\bar{B}=\{ (x,y); x\in B, y\in D(x)\}$,
and where $\bar{T}:\bar{B}\to \bar{B}$ is defined by
$$
\bar{T}(x,y)=\left( T(x), V^{\#}(k(x))(y)\right) .
$$
If $\bar{T}$ is measurable, then $(\bar{B},\bar{T})$ equipped with
the obvious product $\sigma$-algebra is an invertible dynamical
system with invariant density $K$. This dynamical system
$(\bar{B},\bar{T})$ is the natural extension of $(B,T)$.
\end{theorem}

In our setting we have by construction that
$\bar{B}=\Omega_{\alpha}$. So to apply Schweiger's Theorem, we
need to find a backward algorithm $(B^{\#},T^{\#})$, such that
\begin{equation}\label{Vshap-map}
V^{\#}(k(x))(y)=\frac{1}{d(x)\lambda +\varepsilon (x)\,y}.
\end{equation}
{}From our construction of $\Omega_{\alpha}$, and in particular from
the proof of Theorem~\ref{thm:even case nat ext}, we see that a
natural choice of $B^{\#}$ seems\footnote{In a moment we will see
that this choice needs a small modification in case $\alpha \neq
1/2,\, 1/\lambda$.} to be the interval $[0,1]$, with partition
elements
$$
\Delta^{\#}(-1:1)=[H_2,H_{2p-1}]=\left[ \frac{1}{\lambda},
1\right] ,\qquad \Delta^{\#}( 1:1)=[H_1,H_2]=\left[
\frac{1}{\lambda +1},\frac{1}{\lambda}\right] ,
$$
and for $d\geq 2$, $d\neq d_p(\ell_0), d_p(r_0)$,
$$
\Delta^{\#}(-1:d)=\left[ \frac{1}{d\lambda }, \frac{1}{d\lambda
-H_{2p-2}}\right] ,\qquad \Delta^{\#}( 1:d)=\left[
\frac{1}{d\lambda +H_{2p-1}},\frac{1}{d\lambda}\right] .
$$
On these intervals $\Delta^{\#}(\varepsilon : d)$ we define the
map $T^{\#}$ as
$$
T^{\#}(x)=\frac{\varepsilon}{x}-\varepsilon d\lambda ;
$$
An easy calculation shows that $V^{\#}(k(x))$
satisfies~(\ref{Vshap-map}).\smallskip

For $d=d_p(\ell_0)=d_p(r_0)-1$, the rectangle 
$[\ell_0,r_0)\times[\frac1{(d+1)\lambda},\frac1{d\lambda}]$ is 
partitioned as in Figure~\ref{''blow-up''} (see also 
Theorem~\ref{thm:evenorder}, the proof of
Theorem~\ref{thm:even case nat ext}, and Figure~\ref{fig:nat ext q=6}, 
for the case $q=6$, $\alpha=0.53$).

\begin{figure}[h!t]
\includegraphics{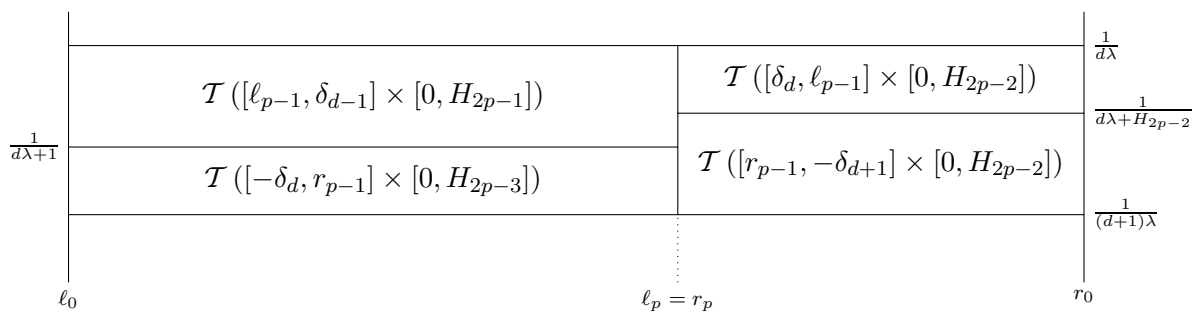}
\caption{``Blow-up'' of the relevant part of $\Omega_\alpha$ for even
$q$, $1/2<\alpha<1/\lambda$.}\label{''blow-up''}
\end{figure}

Note that this case does not occur in case $\alpha =1/2$ or $\alpha
=1/\lambda$; see also Examples~21.3 in~\cite{[Schw]}. Recall, that
$H_{2p-1}=1$, $H_{2p-2}=\lambda /2$, and $H_{2p-3}=\lambda -1$; see
Theorem~\ref{thm:even case nat ext}.\smallskip\

We are left to define $T^{\#}$ on the interval $\left[
\frac{1}{(d+1)\lambda},\frac{1}{d\lambda}\right]$. Clearly, on the
interval $\left[ \frac{1}{(d+1)\lambda},\frac{1}{d\lambda
+1}\right]$ we should define $T^{\#}$ by
$$
T^{\#}(x)=\frac{-1}{x}+(d+1)\lambda ,
$$
while on $\left[ \frac{1}{d\lambda
+H_{2p-2}},\frac{1}{d\lambda}\right]=\left[ \frac{1}{d\lambda
+\lambda /2},\frac{1}{d\lambda}\right]$, the map $T^{\#}$ should
be given by
$$
T^{\#}(x)=\frac{1}{x}-d\lambda .
$$
This leaves us how to define $T^{\#}$ on $\left[ \frac{1}{d\lambda
+1},\frac{1}{d\lambda +H_{2p-2}}\right]$. The problem here is, that
$y$'s in this interval have``two names''; see also
Figure~\ref{''blow-up''}.\smallskip\

\noindent In order to overcome this problem, consider the following
isomorphic copy $(\Omega_{\alpha}^{*},\nu_{\alpha}^{*},{\mathcal
T}^{*})$ of $(\Omega_{\alpha},\nu_{\alpha},{\mathcal T})$; let
$$
\Omega_{\alpha}^{*}=\left( \Omega_{\alpha}\setminus
(\ell_p,r_0]\times \left[\frac{1}{d\lambda +1}, \frac{1}{d\lambda
+\lambda /2} \right] \right) \cup \left( [\ell_p,r_0]\times \left[
\frac{-1}{d\lambda +\lambda /2}, \frac{-1}{d\lambda
+1}\right]\right) ,
$$
let ${\mathcal T}^{*}: \Omega_{\alpha}^{*}\to \Omega_{\alpha}^{*}$
be defined by
$$
{\mathcal T}^{*}(x,y)=s\left( {\mathcal T}(s^{-1}(x,y))\right)
,\qquad (x,y)\in \Omega_{\alpha}^{*},
$$
where
$$
s(x,y)=\begin{cases} (x,y) & \text{if } (x,y)\in
\Omega_{\alpha}\setminus \left( (\ell_p,r_0]\times \left[
\frac{1}{d\lambda +1}, \frac{1}{d\lambda +\lambda /2}\right]
\right) \\
(x,-y) & \text{if } (x,y)\in (\ell_p,r_0]\times \left[
\frac{1}{d\lambda +1}, \frac{1}{d\lambda +\lambda /2} \right],
\end{cases}
$$
and set
$$
\nu_{\alpha}^{*}(A)=\nu_{\alpha}(s^{-1}(A)),\quad \text{for every
Lebesgue set } A\subset \Omega_{\alpha}^*.
$$

Setting
$$
B^{\#}=[0,1]\cup \left( (\ell_p,r_0]\times \left[ \frac{-1}{d\lambda
+\lambda /2}, \frac{-1}{d\lambda +1} \right]\right) ,
$$
and  defining $T^{\#}$ on $\left[ \frac{1}{d\lambda +1},
\frac{1}{d\lambda +\lambda /2}\right]$ by
$$
T^{\#}(x)=\frac{1}{x}-d\lambda ,
$$
and $T^{\#}$ on $\left[ \frac{-1}{d\lambda +\lambda /2},
\frac{-1}{d\lambda+1}\right]$ by
$$
T^{\#}(x)=\begin{cases} \frac{1}{x}+(d+1)\lambda  & \text{if }
\frac{1}{x}+(d+1)\lambda \not\in \left[ \frac{1}{d\lambda +1},
\frac{1}{d\lambda +\lambda /2}\right] \\
\frac{-1}{x}-(d+1)\lambda  & \text{if } \frac{1}{x}+(d+1)\lambda \in
\left[ \frac{1}{d\lambda +1}, \frac{1}{d\lambda +\lambda /2}\right]
,
\end{cases}
$$
we find that $(B^{\#},T^{\#})$ is the dual fibred system with
respect to $(I_{q,\lambda},T_{\alpha})$. We already saw that the
density $K$ of the invariant measure is given by
$$
\frac{C_{q,\alpha}}{(1+xy)^2},\quad \text{with normalizing
constant } C_{q,\alpha}=1/\log\frac{1+\cos\frac\pi q}{\sin\frac\pi
q};
$$
see Proposition~\ref{lem:norm const even}. Thus it follows from
Schweiger's theorem that $(\Omega_{\alpha}, \bar{\mathcal B},
\nu_{\alpha}, {\mathcal T}_{\alpha})$ is the natural extension of
$(I_{q,\alpha}, {\mathcal B}, \mu_{\alpha},T_{\alpha})$.
\end{proof}

\subsection{Odd indices}
Let $q=2h+3$, $h\in\mathbb N$. The $\ell_n$'s and $r_n$'s are
ordered in the following way.
\begin{theorem}\label{thm:oddorder}
Let $q=2h+3$, $h\in\mathbb N$, let the sequences $(\ell_n)_{n\geq
0}$ and $(r_n)_{n\geq 0}$ be defined as above, and let
$\rho=\frac{\lambda-2+\sqrt{\lambda^2-4\lambda+8}}2$. Then we have
the following cases:
$$
\begin{array}{cl}
\alpha=1/2: &
\ell_0<r_{h+1}=\ell_{h+1}<r_1=\ell_1<\cdots<r_{2h-1}=\ell_{2h-1}<
r_{h-1}=\ell_{h-1} \\ & \ <r_{2h}=\ell_{2h}<-\delta_1<r_h=\ell_h<
-\delta_2<r_{2h+1}=\ell_{2h+1}=0<r_0 \vspace{2mm}\\
1/2<\alpha<\rho/\lambda: &
\ell_0<r_{h+1}<\ell_{h+1}<r_1<\cdots<\ell_{h-2}<r_{2h-1}<\ell_{2h-1}<
r_{h-1} \\ & \ <\ell_{h-1}<r_{2h}<\ell_{2h}<-\delta_1<r_h<\ell_h<
-\delta_2<r_{2h+1}<0<\ell_{2h+1}<r_0 \\
& \text{Furthermore, we have $\ell_{2h+2}=r_{2h+2}$ and
$d_{2h+2}(r_0)=d_{2h+2}(\ell_0)+1$.} \vspace{2mm}\\
\alpha=\rho/\lambda: &
\ell_0=r_{h+1}<\ell_{h+1}=r_1<\cdots<\ell_{h-1}<r_h=-\delta_1<\ell_h<
-\delta_2<0<r_0 \vspace{2mm}\\
\rho/\lambda<\alpha<1/\lambda: &
\ell_0<r_1<\cdots<\ell_{h-1}<r_h<-\delta_1<\ell_h<0<r_{h+1}<r_0 \\
& \text{Furthermore, we have $\ell_{h+1}=r_{h+2}$ and
$d_{h+1}(\ell_0)=d_{h+2}(r_0)+1$.} \vspace{2mm}\\
\alpha=1/\lambda: &
\ell_0=r_1<\cdots<\ell_{h-1}=r_h<-\delta_1<\ell_h=r_{h+1}=0<r_0
\end{array}
$$
\end{theorem}

\begin{proof}
In \cite{[BKS]}, Section~3.2 (see also the introduction of this
section), it was shown that
$$
\phi_0=-\frac\lambda2<\phi_1<\cdots<\phi_{h-2}<\phi_{h-1}<
-\frac2{3\lambda}<\phi_h<-\frac2{5\lambda},
$$
$\phi_0<\phi_{h+1}<\phi_1$, and
$$
\phi_{h+1}<\phi_{h+2}<\cdots<\phi_{2h}=-\frac1\lambda<-\frac2{3\lambda}
<\phi_{2h+1}=0.
$$
In view of this and Lemma~\ref{lemma1}, we therefore have that
$\phi_{h-1}\le\ell_{h-1}\le r_h<\phi_h$. An important question is
to know, where $-\delta_1$ is located. Since $3/2\le1+\alpha$, we
have $\phi_{h-1}<-\delta_1$. For $q=2h+3$, we have
$\sin\frac{(h+1)\pi}q=\sin\frac{(h+2)\pi}q$, thus
$$
B_{h+1}=B_{h+2},\qquad B_h=(\lambda-1)B_{h+1},\qquad
B_{h-1}=(\lambda^2-\lambda-1)B_{h+1}.
$$
Hence we obtain, by (\ref{ln}),
$$
\ell_{h-1}=-\frac{\alpha\lambda B_h-B_{h+1}}{\alpha\lambda
B_{h-1}-B_h} =-\frac{1-\alpha\lambda(\lambda-1)}
{\lambda-1-\alpha\lambda(\lambda^2-\lambda-1)} < -\delta_1.
$$
The position of $r_h$ with respect to $-\delta_1$ leads us to
distinguish between the possible cases. We have that
$$
r_h = -\frac{B_{h+1}\alpha\lambda-B_h}{B_h\alpha\lambda-B_{h-1}} =
-\frac{1-(1-\alpha)\lambda}{1-(1-\alpha)\lambda(\lambda-1)} <
-\delta_1
$$
if and only if $\alpha^2\lambda^2+\alpha\lambda(2-\lambda)-1>0$,
i.e.,
$\alpha\lambda>\frac{\lambda-2+\sqrt{\lambda^2-4\lambda+8}}2=\rho$.
Note that
$$
\frac12<\frac{\lambda-2+\sqrt{\lambda^2-4\lambda+8}}{2\lambda}<
\frac1\lambda\qquad\text{for }0<\lambda<2.
$$

Assume first that $\alpha>\rho/\lambda$. Then we have that
$r_h<-\delta_1$, from which it immediately follows that
\begin{gather*}
r_{h+1} =
\frac{B_{h+1}-B_{h+2}\alpha\lambda}{B_{h+1}\alpha\lambda-B_h}
= \frac{1-\alpha\lambda}{1-(1-\alpha)\lambda} \ge 0, \\
\ell_h =
\frac{B_{h+1}\alpha\lambda-B_{h+2}}{B_{h+1}-B_h\alpha\lambda} =
-\frac{1-\alpha\lambda}{1-\alpha\lambda(\lambda-1)} \le 0.
\end{gather*}
If $\alpha<1/\lambda$, then $|1/\ell_h|-|1/r_{h+1}|=\lambda$,
hence $\ell_{h+1}=r_{h+2}$ and $d_1(r_{h+1})=d_1(\ell_{h})+1$. In
case $\alpha=1/\lambda$, then $r_{h+1}=\ell_h=0$. Hence the last
two cases are proved.

It remains to consider $1/2<\alpha\le\rho/\lambda$. Now we have
that $-\delta_1\le r_h(<\phi_h)$, with $-\delta_1=r_h$ if and only
if $\alpha=\rho/\lambda$. Consequently, we immediately find that
$$
r_0=[\, +1:1,(-1:1)^{h-1},-1:2,(-1:1)^h,\ldots\,]\,.
$$
To see that $\ell_h<-\delta_2$, note that this is equivalent to
$$
\alpha^2\lambda^2-\alpha\lambda^2+2\alpha\lambda-1 <
2(\lambda-1)(1-\alpha\lambda),
$$
which holds because of the assumption
$\alpha^2\lambda^2+\alpha\lambda(2-\lambda)-1\le0$. This
assumption also implies that $\ell_{h+1}\le r_1$, where again
equality holds if and only if $\alpha=\rho/\lambda$. This proves
the case $\alpha=\rho/\lambda$.

For $\alpha<\rho/\lambda$, we have
$$
\ell_0=[(-1:1)^h,-1:2,(-1:1)^h,\ldots\,]\,,
$$
hence the convergents of $\ell_0$ satisfy
$-R_{h-1}=B_{h-1}=(\lambda^2-\lambda-1)B_{h+1}$,
$-R_h=(\lambda-1)B_{h+1}$ and
$$
R_{h+1}=2\lambda R_h-R_{h-1}=-(\lambda^2-\lambda+1)B_{h+1}.
$$
The recurrence $R_{h+n+1}=\lambda R_{h+n}-R_{h+n-1}$ for
$n=1,2,\ldots,h$, yields
$$
-R_{h+n} = (B_{n+2}-B_{n+1}+2B_n)B_{h+1}\quad \text{for }
n=0,1,\ldots,h+1.
$$
For the $S_n$'s, we have similarly that
$S_{h-1}=(\lambda-1)B_{h+1}$, $S_h=B_{h+1}$, thus
$S_{h+1}=2\lambda S_h-S_{h-1}=(\lambda+1)B_{h+1}$ and
$$
S_{h+n} = (B_{n+1}+B_n)B_{h+1}\qquad \text{for }n=0,1,\ldots,h+1.
$$
By (\ref{xvsTn}), we obtain, for $n=0,1,\ldots,h+1$,
\begin{equation}\label{lnh}
\ell_{h+n}=-\frac{(B_{n+1}+B_n)(\alpha-1)\lambda+B_{n+2}-B_{n+1}+2B_n}
{(B_n+B_{n-1})(\alpha-1)\lambda+B_{n+1}-B_n+2B_{n-1}}\,.
\end{equation}
For the convergents of $r_0$, only the sign of the $R_n$ is
different and we get
\begin{equation}\label{rnh}
r_{h+n} =
-\frac{(B_{n+1}+B_n)\alpha\lambda-(B_{n+2}-B_{n+1}+2B_n)}
{(B_n+B_{n-1})\alpha\lambda-(B_{n+1}-B_n+2B_{n-1})}\,.
\end{equation}
This yields that
$$
r_{2h+1} = -\frac{(2\alpha-1)\lambda}{\alpha\lambda^2-2\lambda+2}\
(<0),\qquad \ell_{2h+1} =
\frac{(2\alpha-1)\lambda}{(1-\alpha)\lambda^2-2\lambda+2}\ (>0),
$$
hence $|1/r_{2h+1}|-|1/\ell_{2h+1}|=\lambda$,
$d_1(r_{2h+1})=d_1(\ell_{2h+1})+1$, and the theorem is proved.
\end{proof}

For the construction of the natural extension, we have to
distinguish between the different cases of the previous theorem.
Consider first $\alpha>\rho/\lambda$.

\begin{theorem}\label{thm:nat ext odd1}
Let $q=2h+3$ with $h\ge1$. Then the system of relations
$$
\left\{\begin{array}{rl}
(\mathcal R_1): & \qquad H_1=1/(\lambda+H_{2h+2}) \\
(\mathcal R_2): & \qquad H_2=1/\lambda \\
(\mathcal R_n): & \qquad H_n=1/(\lambda-H_{n-2})\quad \text{for }
n=3,4,\ldots,2h+2 \\
(\mathcal R_{2h+3}): & \qquad H_{2h+1}=\lambda/2 \\
(\mathcal R_{2h+4}): & \qquad H_{2h}+H_{2h+2}=\lambda
\end{array}\right.
$$
admits the (unique) solution
\begin{gather*}
H_{2n} = -\phi_{2h+1-n} = \frac{B_n}{B_{n+1}} =
\frac{\sin\frac{n\pi}q}{\sin\frac{(n+1)\pi}q}\ \text{ for }
n=1,2,\ldots,h+1, \\
H_{2n-1} = \frac{B_{n-1}+B_n}{B_n+B_{n+1}} =
\frac{\sin\frac{(n-1)\pi}q+\sin\frac{n\pi}q}
{\sin\frac{n\pi}q+\sin\frac{(n+1)\pi}q} \ \text{ for }
n=1,2,\ldots,h+1,
\end{gather*}
in particular $H_{2h+2}=1$.

Let $\rho/\lambda<\alpha\le1/\lambda$ and
$\Omega_\alpha=\bigcup_{n=1}^{2h+2}J_n\times[0,H_n]$ with
$J_{2n-1}=[\ell_{n-1},r_n)$, $J_{2n}=[r_n,\ell_n)$ for
$n=1,2,\ldots,h$, $J_{2h+1}=[\ell_h,r_{h+1})$ and
$J_{2h+2}=[r_{h+1},r_0)$. Then the map $\mathcal T_\alpha:\,
\Omega_\alpha \rightarrow \Omega_\alpha$ given by (\ref{natural
extension map}) is bijective off of a set of Lebesgue measure
zero.
\end{theorem}

\begin{remark0}
The case $q=3$, $\rho/\lambda\leq\alpha\leq 1/\lambda$, which is
the case of Nakada's $\alpha$-expansions for
$(\sqrt{5}-1)/2\leq\alpha\leq 1$, has been dealt with
in~\cite{[N1]}; see also~\cite{[NIT]}, \cite{[TI]},
and~\cite{[K1],[K2]}.
\end{remark0}

The proof of Theorem~\ref{thm:nat ext odd1} is very similar to
that of Theorem~\ref{thm:even case nat ext} and therefore omitted,
see also Figure~\ref{fig:nat ext q=5 case1}. In case
$\alpha=1/\lambda$, the intervals $J_{2n-1}$ are empty.\medskip\

\begin{figure}[h!t]
\includegraphics{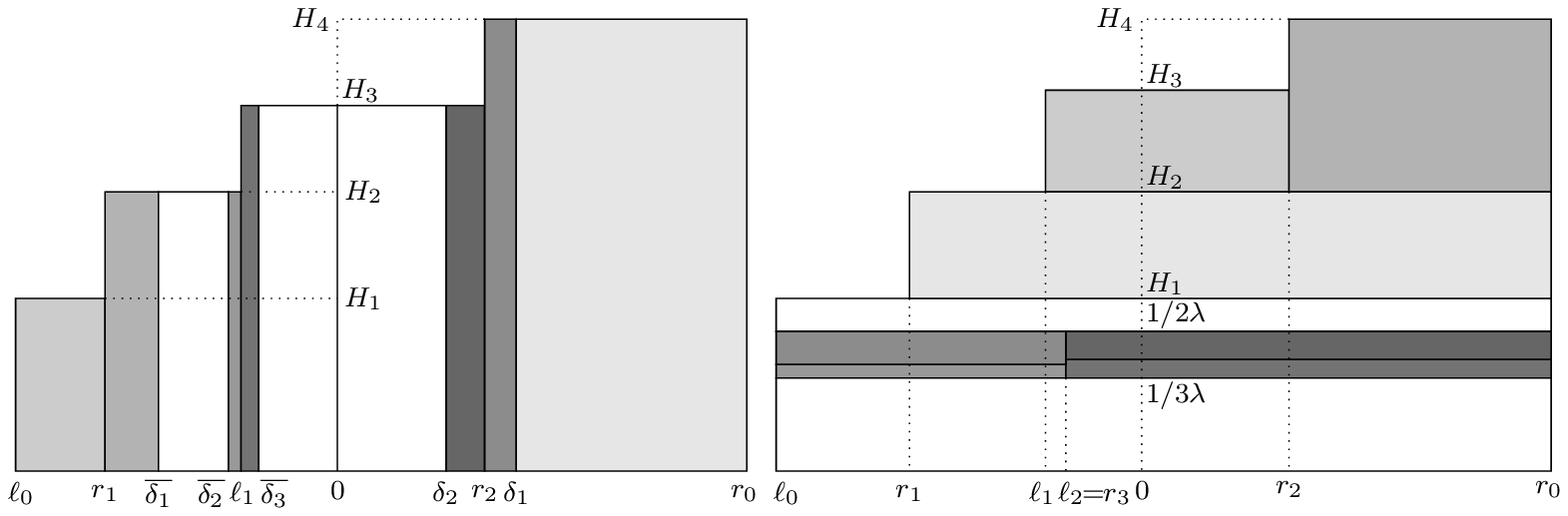}
\caption[natural extension q=5, first case]{The natural extension
domain $\Omega_\alpha$ (left) and its image under $\mathcal
T_\alpha$ (right) of the $\alpha$-Rosen continued fraction
($\overline{\delta_n}=-\delta_n$); here $q=5$, $\alpha=0.56$,
$d_{h+1}(\ell_0)=3$, $d_{h+2}(r_0)=2$.} \label{fig:nat ext q=5
case1}
\end{figure}

Once more, a Jacobian calculation shows that $\mathcal T_\alpha$
preserves the probability measure $\nu_\alpha$ with density
$$
\frac{C_{q,\alpha}}{(1+xy)^2} \, ,
$$
where $C_{q,\alpha}$ is a normalizing constant given by the
following proposition.

\begin{proposition}\label{lem:norm const odd1}
If $q=2h+3$ and $\rho/\lambda<\alpha\le1/\lambda$, then the
normalizing constant is
$$
C_{q,\alpha}=1/\log\frac{1+\alpha\lambda}{\sqrt{2-\lambda}} =
1/\log\frac{1+2\alpha\cos\frac\pi q}{2\sin\frac\pi{2q}}\,.
$$
\end{proposition}

\begin{proof}
Integration gives
$$
C_{q,\alpha} = 1/\log \left( \frac{1+r_0}{1+r_{h+1}}
\prod_{n=1}^{h+1} \frac{1+r_nH_{2n-1}}{1+\ell_{n-1}H_{2n-1}}
\prod_{n=1}^h \frac{1+\ell_nH_{2n}}{1+r_nH_{2n}} \right)\,.
$$
Using (\ref{ln}) and (\ref{rn}), we find
$$
\frac{1+r_nH_{2n-1}}{1+\ell_{n-1}H_{2n-1}} =
\frac{B_n-B_{n-1}\alpha\lambda}{B_n\alpha\lambda-B_{n-1}}\,,\qquad
\frac{1+\ell_nH_{2n}}{1+r_nH_{2n}} =
\frac{B_n\alpha\lambda-B_{n-1}}{B_{n+1}-B_n\alpha\lambda}\,,
$$
and
$$
\frac{1+r_0}{1+r_{h+1}} =
\frac{(1+\alpha\lambda)(B_{h+1}\alpha\lambda-B_h)}{B_{h+1}-B_h} =
(1+\alpha\lambda)(B_{h+1}\alpha\lambda-B_h)B_{h+1}\,,
$$
where we have used
$$
(2-\lambda)B_{h+1}^2 = \frac{2(1-\cos\frac\pi
q)\sin^2\frac{(h+1)\pi}q}{\sin^2\frac\pi q} =
\frac{4\sin^2\frac\pi{2q}\cos^2\frac\pi{2q}}
{4\sin^2\frac\pi{2q}\cos^2\frac\pi{2q}} = 1.
$$
Putting everything together, we obtain
$$
C_{q,\alpha} = 1/\log ((1+\alpha\lambda)B_{h+1}) = 1/\log
\frac{1+\alpha\lambda}{\sqrt{2-\lambda}} = 1/\log
\frac{1+\alpha\lambda}{2\sin\frac\pi{2q}}\,.
$$
\end{proof}

\begin{remark0}
Note that for $q=3$ this result confirms Nakada's result
from~\cite{[N1]} for $\alpha$ between $(\sqrt{5}-1)/2$ and $1$; in
this case, the normalizing constant is indeed $1/\log(1+\alpha)$.
\end{remark0}

Now consider $\alpha<\rho/\lambda$.

\begin{theorem}\label{thm:nat ext odd2}
Let $q=2h+3$ with $h\ge1$. Then the system of relations
$$
\left\{\begin{array}{rl}
(\mathcal R_1): & \qquad H_1=1/(2\lambda-H_{4h-1}) \\
(\mathcal R_2): & \qquad H_2=1/(2\lambda-H_{4h}) \\
(\mathcal R_3): & \qquad H_3=1/(\lambda+H_{4h+3}) \\
(\mathcal R_4): & \qquad H_4=1/\lambda \\
(\mathcal R_n): & \qquad H_n=1/(\lambda-H_{n-4})\qquad \text{for }
n=5,6,\dots,4h+3 \\
(\mathcal R_{4h+4}): & \qquad H_{4h+2}=\lambda/2 \\
(\mathcal R_{4h+5}): & \qquad H_{4h+1}+H_{4h+3}=\lambda
\end{array}\right.
$$
admits the (unique) solution
\begin{gather*}
H_{4n}=\frac{B_n}{B_{n+1}},\qquad
H_{4n-2}=\frac{B_{n-1}+B_n}{B_n+B_{n+1}}, \\
H_{4h+3-4n}=\frac{B_{n+1}\rho-B_n}{B_n\rho-B_{n-1}}, \qquad
H_{4h+1-4n}=\frac{B_{n+1}\rho-B_{n+2}}{B_n\rho-B_{n+1}},
\end{gather*}
in particular $H_{4h+3}=\rho$.

Let $1/2\le\alpha<\rho/\lambda$, and
$\Omega_\alpha=\bigcup_{n=1}^{4h+3}J_n\times [0,H_n]$ with
$J_{4n-3}=[\ell_{n-1},r_{h+n})$, $J_{4n-2}=[r_{h+n},\ell_{h+n})$,
$J_{4n-1}=[\ell_{h+n},r_n)$, $J_{4n}=[r_n,\ell_n)$ for
$n=1,2,\ldots,h$, $J_{4h+1}=[\ell_h,r_{2h+1})$,
$J_{4h+2}=[r_{2h+1},\ell_{2h+1})$ and
$J_{4h+3}=[\ell_{2h+1},r_0)$. Then the map $\mathcal T_\alpha:\,
\Omega_\alpha \rightarrow \Omega_\alpha$ given by (\ref{natural
extension map}) is bijective off of a set of Lebesgue measure
zero.
\end{theorem}

\begin{proof}
The proof of the bijectivity runs along the same lines as the
proof of Theorem~\ref{thm:even case nat ext} and is therefore
omitted, see also Figure~\ref{fig:nat ext q=5 case2}.

The $H_{4n}$'s are determined by $(\mathcal R_4),(\mathcal
R_8),\ldots,(\mathcal R_{4h})$. The $H_{4n-2}$'s are determined by
$(\mathcal R_2),(\mathcal R_6),\ldots,(\mathcal R_{4h+2})$ and
$(\mathcal R_{4h+4})$. By $(\mathcal R_3),(\mathcal
R_7),\ldots,(\mathcal R_{4h+3})$, we obtain
$H_{4h+3-4n}=\frac{B_{n+1}H_{4h+3}-B_n}{B_nH_{4h+3}-B_{n-1}}$ and
$$
\frac1{\lambda+H_{4h+3}} = H_3 =
\frac{B_{h+1}H_{4h+3}-B_h}{B_hH_{4h+3}-B_{h-1}} = \frac{H_{4h+3}-
(\lambda-1)}{(\lambda-1)H_{4h+3}-(\lambda^2-\lambda-1)}\,,
$$
thus $H_{4h+3}^2+(2-\lambda)H_{4h+3}-1=0$, i.e. $H_{4h+3}=\rho$.
Finally, the $H_{4h+1-4n}$'s are determined by $(\mathcal
R_1),(\mathcal R_5),\ldots,(\mathcal R_{4h+5})$. For $\alpha=1/2$,
the intervals $J_{4n}$ and $J_{4n-2}$ are empty.
\end{proof}

\begin{figure}[h!t]
\includegraphics{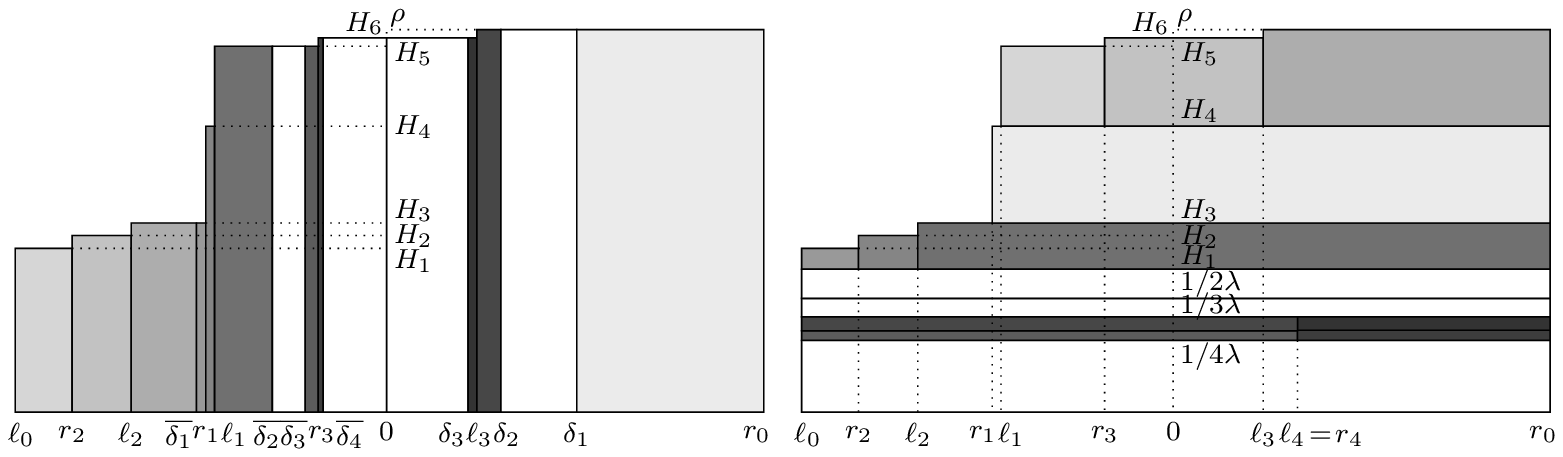}
\caption[natural extension q=5, second case]{The natural extension
domain $\Omega_\alpha$ (left) and its image under $\mathcal
T_\alpha$ (right) of the $\alpha$-Rosen continued fraction
($\overline{\delta_n}=-\delta_n$); here $q=5$, $\alpha=0.5038$,
$d_{2h+2}(\ell_0)=2$, $d_{2h+2}(r_0)=3$.} \label{fig:nat ext q=5
case2}
\end{figure}

Again, $\mathcal T_\alpha$ preserves the probability measure
$\nu_\alpha$ with density $C_{q,\alpha}/(1+xy)^2$, where
$C_{q,\alpha}$ is a normalizing constant given by the following
proposition.
\begin{proposition}\label{lem:norm const odd2}
If $q=2h\!+\!3$ and $1/2\le\alpha<\rho/\lambda$, then the
normalizing constant is
$$
C_{q,\alpha} = 1/\log\frac{1+\rho}{\sqrt{2-\lambda}} =
1/\log\frac{1+\rho}{2\sin\frac\pi{2q}}.
$$
\end{proposition}

\begin{proof}
Integration yields that $C_{q,\alpha}$ is equal to
$$
1/\log \left( \frac{1+r_0\rho}{1+\ell_{2h+1}\rho}
\prod_{j=1}^{h+1} \frac{1+r_{h+j}H_{4j-3}}{1+\ell_{j-1}H_{4j-3}}
\frac{1+\ell_{h+j}H_{4j-2}}{1+r_{h+j}H_{4j-2}} \prod_{j=1}^h
\frac{1+r_jH_{4j-1}}{1+\ell_{h+j}H_{4j-1}}
\frac{1+\ell_jH_{4j}}{1+r_jH_{4j}} \right) .
$$
Using (\ref{ln}), (\ref{rn}), (\ref{lnh}), (\ref{rnh}) and
Lemma~\ref{lem:B-identity}, we find
\begin{gather*}
\frac{1+r_{h+j}H_{4j-3}}{1+\ell_{j-1}H_{4j-3}} =
\frac{(\lambda(1-2\alpha)\rho+\alpha\lambda^2-\lambda^2+2\lambda-2)
(B_j-B_{j-1}\alpha\lambda)}{((\alpha\lambda-1)\rho-\alpha\lambda+
\lambda-1)((B_j+B_{j-1})\alpha\lambda-B_{j+1}+B_j-2B_{j-1})} \\
\frac{1+\ell_{h+j}H_{4j-2}}{1+r_{h+j}H_{4j-2}} =
\frac{(B_j+B_{j-1})\alpha\lambda-B_{j+1}+B_j-2B_{j-1}}
{2B_j-B_{j-1}+B_{j-2}-(B_j+B_{j-1})\alpha\lambda} \\
\frac{1+r_jH_{4j-1}}{1+\ell_{h+j}H_{4j-1}} =
\frac{((\alpha\lambda-1)\rho-\alpha\lambda+\lambda-1)
(2B_j-B_{j-1}+B_{j-2}-(B_j+B_{j-1})\alpha\lambda)}
{(\lambda(1-2\alpha)\rho+\alpha\lambda^2-\lambda^2+2\lambda-1)
(B_j\alpha\lambda-B_{j-1}) }\\
\frac{1+\ell_jH_{4j}}{1+r_jH_{4j}} =
\frac{B_j\alpha\lambda-B_{j-1}}{B_{j+1}-B_j\alpha\lambda}
\end{gather*}
and
$$
\frac{1+r_0\rho}{1+\ell_{2h+1}\rho} = \frac{(1+\alpha\lambda\rho)
((2B_{h+1}-B_h+B_{h-1}-(B_{h+1}+B_h)\alpha\lambda)}
{((2\alpha-1)\lambda\rho-\alpha\lambda^2+\lambda^2+2\lambda-2)B_{h+1}}\,.
$$
Putting everything together, we obtain that
$$
C_{q,\alpha} = 1/\log
\frac{\left(1+\alpha\lambda\rho\right)\sqrt{2-\lambda}}
{1+\alpha\lambda-\lambda+(1-\alpha\lambda)\rho}=
1/\log\frac{1+\rho}{\sqrt{2-\lambda}}\,.
$$
\end{proof}

This confirms again Nakada's result for $q=3$, i.e.,
$C_{3,\alpha}=1/\log\frac{\sqrt5+1}2$ for
$\frac12\le\alpha<\frac{\sqrt5-1}2$.

\smallskip
The case $\alpha=\rho/\lambda$ is slightly different from both
other cases (similarly to $\alpha=1/\lambda$ for even $q$).

\begin{theorem}
Let $q=2h+3$ with $h\ge1$, $\alpha=\rho/\lambda$ and
$$
\Omega_\rho = \bigcup_{j=1}^h \left([\ell_{j-1},r_j)\times
\left[0,\frac{B_{j-1}+B_j}{B_j+B_{j+1}}\right]\right) \cup
\left([r_j,\ell_j)\times\left[0,\frac{B_j}{B_{j+1}}\right]\right)
\cup \left([\ell_h,\rho)\times\left[0,\frac\lambda2\right]\right).
$$
Then $\mathcal T_\rho:\, \Omega_\rho\rightarrow\Omega_\rho$ is
bijective off of a set of Lebesgue measure zero.
\end{theorem}

The normalizing constant in this case is $C_{q,\rho/\lambda} =
1/\log\frac{1+\rho}{\sqrt{2-\lambda}}$ as above. As in the even
case, we set $\mu_{\alpha}$ the projection of $\nu_{\alpha}$ on the
first coordinate, $\bar{\mathcal B}$ be the restriction of the
two-dimensional $\sigma$-algebra on $\Omega_{\alpha}$, and
${\mathcal B}$ be the Lebesgue $\sigma$-algebra on
$I_{q,\alpha}=[\lambda (\alpha -1),\alpha \lambda ]$. We have the
following theorem, whose proof is similar to the proof of
Theorem~ref{thm:natural-extension-even-case} in the even case.
\begin{theorem}\label{thm:natural-extension-odd-case}
Let $q\geq 3$, $q=2h+1$, and let $\frac{1}{2}\leq \alpha \leq
\frac{1}{\lambda}$. Then the dynamical system $(\Omega_{\alpha},
\bar{\mathcal B}, \nu_{\alpha},{\mathcal T}_{\alpha})$ is the
natural extension of the dynamical system $(I_{q,\alpha}, {\mathcal
B}, \mu_{\alpha},T_{\alpha})$.
\end{theorem}

\begin{figure}[h!t]
\includegraphics{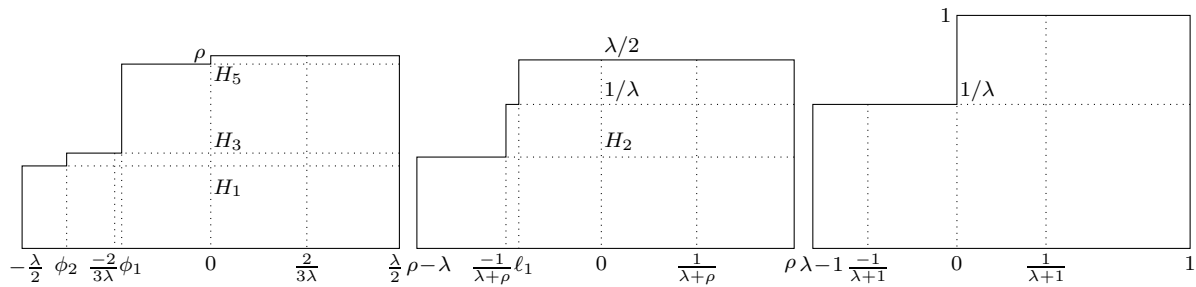}
\caption[1/2, rho/lambda and 1/lambda]{$\Omega_{1/2}$ (left),
$\Omega_{\rho/\lambda}$ (middle) and $\Omega_{1/\lambda}$ (right);
here $q=5$.} \label{fig:nat ext 1/2,rho/lambda,1/lambda}
\end{figure}

\subsection{Convergence of the continued fractions}\label{conv}
Now we can prove easily that the $\alpha$-Rosen continued
fractions converge. If $T_\alpha^n(x)=0$ for some $n\ge0$, then
this is clear. Therefore assume that $T_\alpha^n(x)\ne 0$ for all
$n\ge0$. Setting $(t_n,v_n)=\mathcal{T}_{\alpha}^n(x,0)$, it
follows directly from the definition (\ref{natural extension map})
of $\mathcal{T}_{\alpha}$ that
$v_n=[\,1:d_n,\,\varepsilon_n:d_{n-1},\dots\,,\varepsilon_2:d_1\,]$.
Furthermore, an immediate consequence of~(\ref{recurrence
relations}) is that
$S_{n-1}/S_n=[\,1:d_n,\,\varepsilon_n:d_{n-1},\dots\,,
\varepsilon_2:d_1\,]$, i.e., $v_n=S_{n-1}/S_n$.

Theorems~\ref{thm:even case nat ext}, \ref{thm:nat ext odd1} and
\ref{thm:nat ext odd2} (see also Figures~\ref{fig:nat ext q=6},
\ref{fig:nat ext q=5 case1} and \ref{fig:nat ext q=5 case2}) show
that $v_n\le1$, i.e., $S_n\ge S_{n-1}$, and that $S_n=S_{n-1}$ if
and only if $q=2h+3$, $n=h+1$, $d_1=1$,
$(\varepsilon_i:d_i)=(-1:1)$ for $i=2,3,\ldots,h+1$, which is
possible only if $\alpha>\rho/\lambda$. Furthermore we have that
$v_{n-1}v_n\le 1/c$ for some constant $c>1$, i.e., $S_n\ge
cS_{n-2}$. It follows from (\ref{xvsTn}) that
\begin{equation}\label{thetan}
\left|x-\frac{R_n}{S_n}\right|=
\left|\frac{T_\alpha^n(x)(R_{n-1}S_n-R_nS_{n-1})}
{S_n(S_n+T_\alpha^n(x)S_{n-1})}\right|=\frac{|t_n|}{S_n^2(1+t_nv_n)}\le
\frac{\alpha\lambda}{(1+\alpha\lambda-\lambda)S_n^2}\,,
\end{equation}
hence the $\alpha$-Rosen convergents $R_n/S_n$ converge to $x$ as
$n\to\infty$.


\section{Mixing properties of $\alpha$-Rosen fractions}
\label{mixing properties} In case $q$ is even, and
$\alpha=1/\lambda$, we saw in the previous section that there is a
simple relation between $\Omega_{1/2}$ and $\Omega_{1/\lambda}$;
see also Figure~\ref{fig:nat ext 1/2,1/lambda}. Define in this
case the map $\mathcal M: \Omega_{1/\lambda}\to \Omega_{1/2}$ by
$$
\mathcal M(x,y)=\begin{cases}
(-y,-x) & \text{if $(x,y)\in \Omega_{1/\lambda},\,\, x<0$},\\
(y,x) & \text{if $(x,y)\in \Omega_{1/\lambda},\,\, x\geq 0$}.
\end{cases}
$$
Clearly, $\mathcal M: \Omega_{1/\lambda}\to \Omega_{1/2}$ is
bijective and bi-measurable transformation, and
$\nu_{1/\lambda}(\mathcal M^{-1}(A))=\nu_{1/2}(A)$, for every
Borel set $A\subset\Omega_{1/2}$. By comparing the partitions of
$\mathcal T_{1/\lambda}$ (on $\Omega_{1/\lambda}$) and that of
$\mathcal T_{1/2}^{-1}$ (on $\Omega_{1/2}$), we find that
$$
\mathcal T_{1/\lambda}(x,y)= \mathcal M^{-1}\left(\mathcal
T_{1/2}^{-1} \left(\mathcal M(x,y)\right) \right) ,\qquad (x,y)\in
\Omega_{1/\lambda},\,\, x\neq 0.
$$
This implies that the dynamical systems $(\Omega_{1/2},\bar{\mathcal
B},\nu_{1/2},\mathcal T_{1/2}^{-1})$ and
$(\Omega_{1/\lambda},\bar{\mathcal B},\nu_{1/\lambda},\mathcal
T_{1/\lambda})$ are isomorphic. In~\cite{[BKS]} it was shown that
the dynamical system system $(\Omega_{1/2},\bar{\mathcal
B},\nu_{1/2},\mathcal T_{1/2})$ is weakly Bernoulli with respect to
the natural partition, hence $(\Omega_{1/2},\bar{\mathcal
B},\nu_{1/2},\mathcal T_{1/2})$ and $(\Omega_{1/2},\bar{\mathcal
B},\nu_{1/2},\mathcal T_{1/2}^{-1})$ are isomorphic. As a
consequence we find that the dynamical systems
$(\Omega_{1/2},\bar{\mathcal B},\nu_{1/2},\mathcal T_{1/2})$ and
$(\Omega_{1/\lambda},\bar{\mathcal B},\nu_{1/\lambda},\mathcal
T_{1/\lambda})$ are isomorphic.\medskip\

In this section, we will show that this result also holds for all
$q$ and all $\alpha$ strictly between $1/2$ and $1/\lambda$, using
a result by M.~Rychlik~\cite{[Ry]}. For completeness, we state
explicitly the hypothesis needed for Rychlik's result (the reader
is referred to~\cite{[Ry]} for more details).

Let $X$ be a totally ordered order-complete set. Open intervals
constitute a base of a complete topology in $X$, making $X$ into a
topological space. If $X$ is separable, then $X$ is homeomorphic
with a closed subset of an interval. Let $\mathcal{B}$ be the Borel
$\sigma$-algebra on $X$, and $m$ a fixed regular, Borel probability
measure on $X$ (in our case $m$ will be the normalized Lebesgue
measure restricted to $X$). Let $U\subset X$ be an open dense subset
of $X$, such that $m(U)=1$. Let $S=X\setminus U$, clearly $m(S)=0$.

Let $T:U\to X$ be a continuous map, and $\beta$ a countable family
of closed intervals with disjoint interiors, such that $U\subset
\bigcup \beta$. Furthermore, suppose that for any $B\in \beta$ one
has that $B\cap S$ consists only of endpoints of $B$, and that $T$
restricted to $B\cap U$ admits an extension to a homeomorphism of
$B$ with some interval in $X$. Suppose that $T'(x)\ne 0$ for $x\in
U$, and let $g(x)=1/|T'(x)|$ for $x\in U$, $g|_S=0$. Let
$P:L^1(X,m)\to L^1(X,m)$ be the Perron-Frobenius operator of $T$,
$$
Pf(x)=\sum_{y\in T^{-1}x} g(y) f(y).
$$

In~\cite{[Ry]}, it was proved (among many other things) that, if
$\|g\|_{\infty}<1$ and $\mathrm{Var}\,g<\infty$, then there exist
functions $\varphi_1,\varphi_2,\dots,\varphi_s$ of bounded
variation, such that
\begin{itemize}
\item[(i)] $P\varphi_i=\varphi_i$;
\item[(ii)] $\int \varphi_i\, \mathrm{d}m=1$;
\item[(iii)] There exists a measurable partition $C_1,C_2,\dots,C_s$ of
$X$ with $T^{-1}C_i=C_i$ for $i=1,2,\dots,s$;
\item[(iv)] The dynamical system $(C_i,T_i,\nu_i)$, where
$T_i=T|_{C_i}$ and $\nu_i(B)=\int_B\varphi_i\, \mathrm{d}m$ are
exact, and $\nu_i$ is the unique invariant measure for $T_i$,
absolutely continuous with respect to $m|_{C_i}$.
\end{itemize}
Rychlik also showed that if $s=1$, i.e., if $1$ is the only
eigenvalue of $P$ on the unit circle and if there exists only one
$\varphi\in L^1(X,M)$ with $P\varphi=\varphi$ and $m(\varphi)=1$,
($\varphi\ge 0$), then the natural extension of $(X,T,\nu)$ is
isomorphic to a Bernoulli shift.

Returning to our map $T_{\alpha}$, defined on
$X=I_{q,\alpha}=[\lambda(\alpha-1),\alpha\lambda]$, and using the
same notation as above, we let $m$ be normalized Lebesgue measure on
$X$,
$$
S=\{\lambda(\alpha-1)\}\cup \left\{\pm\frac{1}{\lambda(\alpha
+d)}\,;\ d=1,2,\dots\right\}
$$
and $U=X\setminus S$. Note that $T_\alpha:U\to X$ is continuous,
and that the restriction of $T_\alpha$ to each open interval is
homeomorphic to an interval (in fact to $X$ itself, except for the
first and last interval).

We have that $g(x)=1/|T_{\alpha}^{\prime}(x)|=x^2$ on $U$, hence
$\|g\|_{\infty}<1$ (since $\alpha \neq 1/\lambda$), and
$\mathrm{Var}\,g<\infty$. It is easy, but tedious (cf.~\cite{[DK]}
for a proof of the regular case), to see that $T$ is ergodic, hence
$s=1$, and we can apply Rychlik's result to obtain the following
theorem.
\begin{theorem}\label{thmbernoulli}
The natural extension $(\Omega_{\alpha},\nu_{\alpha},\mathcal
T_{\alpha})$ of $(X_{\alpha}, \mu_{\alpha}, T_{\alpha})$ is weakly
Bernoulli. Hence, the natural extension is isomorphic to any
Bernoulli shift with the same entropy.
\end{theorem}

\section{Metrical properties of `regular' Rosen fractions}
\label{Metrical properties} An important reason to introduce and
study the natural extension of the ergodic system underlying any
continued fraction expansion, is that such a natural extension
facilitates the study of the continued fraction expansion at hand;
see e.g.\ \cite{[DK]} and~\cite{[IK]}, Chapter~4. The following
theorem is a consequence of this; see~\cite{[BJW]}, \cite{[DK]},
or~\cite{[IK]}, Chapter~4.
\begin{theorem}\label{distribution-orbit}
Let $q\geq 3$, and let $1/2\leq \alpha \leq 1/\lambda$. For almost
all $G_q$-irrational numbers $x$, the two-dimensional sequence
$\mathcal
T_\alpha(x,0)=\left(\,T_\alpha^n(x),\,S_{n-1}/S_n\,\right)$,
$n\geq 1$, is distributed over $\Omega_{\alpha}$ according to the
density function $g_{\alpha}$, given by
$$
g_{\alpha}(t,v)= \frac{C_{q,\alpha}}{(1+tv)^2}
$$
for $(t,v)\in \Omega_{\alpha}$, and $g_{\alpha}(t,v)=0$ otherwise.
Here $C_{q,\alpha}$ is the normalizing constant of the $\mathcal
T_{\alpha}$-invariant measure $\nu_{\alpha}$.
\end{theorem}

Due to Proposition~\ref{distribution-orbit}, it is possible to
study the distribution of various sequences related to the
$\alpha$-Rosen expansion of almost every $x\in X_{\alpha}$.
Classical examples of these are the frequency of digits, or the
analogs of various classical results by L\'evy and Khintchine.
However, these results can already be obtained from the projection
$(X_{\alpha}, {\mathcal B}_{\alpha}, \mu_{\alpha}, T_{\alpha})$ of
$(\Omega_\alpha,\bar{\mathcal B}_\alpha,\nu_\alpha,\mathcal
T_\alpha)$
--- which is also ergodic --- and the Ergodic Theorem.
For the distribution of the so-called {\em approximation
coefficients}, the natural extension $(\Omega_\alpha,\bar{\mathcal
B}_\alpha,\nu_\alpha,\mathcal T_\alpha)$ is necessary. These
approximation coefficients $\Theta_n=\Theta_n(x)$, are defined by
\begin{equation}\label{regular-thetas1}
\Theta_n=\Theta_n(x):=S_n^2\left|x-\frac{R_n}{S_n}\right|,\quad
n\ge0,
\end{equation}
where $R_n/S_n$ is the $n$th $\alpha$-Rosen convergent, which is
obtained by truncating the $\alpha$-Rosen expansion.

With $(t_n,v_n)={\mathcal T}_{\alpha}^n(x,0)$, it follows from
(\ref{thetan}), that
\begin{equation}\label{regular-thetas2}
\Theta_n=\frac{\varepsilon_{n+1}t_n}{1+t_nv_n},\quad \text{for
$n\geq 1$}.
\end{equation}
Similarly, since $t_n=\varepsilon_n/t_{n-1}-d_n\lambda$ and
$v_n=S_{n-1}/S_n$, it follows from~(\ref{recurrence relations}) that
\begin{equation}\label{regular-thetas3}
\Theta_{n-1}=\frac{v_n}{1+t_nv_n},\quad \text{for $n\geq 1$}.
\end{equation}

In view of~(\ref{regular-thetas2}) and~(\ref{regular-thetas3}), we
define the map
$$
F(t,v)=\left( \frac{v}{1+tv},\frac{t}{1+tv}\right) =:
(\xi,\eta),\qquad \text{for } tv\neq -1.
$$
It is now easy calculation, see e.g.~\cite{[BKS]}, p.\ 1293, that
due to Proposition~\ref{distribution-orbit} one has for almost all
$x\in X_{\alpha}$ that the sequence
$(\Theta_{n-1}(x),\varepsilon_{n+1}\Theta_n(x))_{n\geq 0}$ is
distributed on $F(\Omega_{\alpha})$ according to the density
function $C_{q,\alpha}/\sqrt{1-4\xi\eta}$. Setting
$$
\Gamma_{\alpha}^+=F(\{ (t,v)\in\Omega_{\alpha} \mid t\geq 0\} )
\quad \text{and}\quad
\Gamma_\alpha^-=F(\{(t,v)\in\Omega_\alpha\mid t\le0\}),
$$
we have found the following theorem.

\begin{theorem}\label{twoTheta's}
Let $q\geq 3$, $1/2\leq \alpha\leq 1/\lambda$, and define the
functions $d_{\alpha}^+$ and $d_{\alpha}^-$ by
\begin{equation}\label{density-d}
d_{\alpha}^\pm(\xi ,\eta )= \frac{C_{q,\alpha}}{\sqrt{1\mp
4\xi\eta}} \quad \text{for}\quad (\xi ,\eta) \in
\Gamma_\alpha^\pm,
\end{equation}
and $d_{\alpha}^\pm(\xi ,\eta )=0$ otherwise. Then the sequence
$(\Theta_{n-1}(x),\Theta_n(x))_{n\ge1}$ lies in the interior of
$\Gamma=\Gamma_{\alpha}^+\cup\Gamma_{\alpha}^-$ for all
$G_q$-irrational numbers $x$, and for almost all $x$ this sequence
is distributed according to the density function $d_{\alpha}$,
where
$$
d_\alpha(\xi,\eta)=d_\alpha^+(\xi,\eta)+d_\alpha^-(\xi,\eta)\,.
$$
By this last statement we mean that, for almost all $x$ and for
all $a,b\geq 0$, the limit
$$
\lim_{N \to \infty} \, \frac{1}{N}\, \# \{j\,:\, 1 \le j \le N,
\Theta_{j-1}(x) < a, \, \Theta_j(x) < b\,\}
$$
exists, and equals
$$
\int_{0}^a\!\!\!\int_{0}^b d_{\alpha}(\xi, \eta)\, {\rm d}\xi\,
{\rm d}\eta.
$$
\end{theorem}

Several corollaries can be drawn from Theorem~\ref{twoTheta's},
see e.g.~[K1], where (for $q=3$) for almost all $x$ the
distributions of the sequences $(\Theta_n)_{n\geq 1}$,
$(\Theta_{n-1} +\Theta_n)_{n\geq 1}$,
$(\Theta_{n-1}-\Theta_n)_{n\geq 1}$ were determined.

Here we only mention the following result for even values of $q$,
a result which was previously obtained in \cite{[BKS]} for both
even and odd values of $q$, and $\alpha =1/2$.

\begin{proposition}\label{lenstra}
Let $q\geq 4$ be an even integer, $1/2\leq \alpha\leq 1/\lambda$,
and let
$$
\mathcal L_{\alpha} := \min \left\{
\frac{\lambda}{\lambda+2},\frac{\lambda
(2-\alpha\lambda^2)}{4-\lambda^2}\right\}.
$$
Then for almost all $G_q$-irrational numbers $x$ and all $c\geq
1/\mathcal L_{\alpha}$, we have that
$$
\lim_{N\to \infty} \frac{1}{N} \# \left\{ n:1\leq n \leq N,\,
\Theta_n(x) < \frac{1}{c} \right\}  = \frac{\lambda
C_{q,\alpha}}{c}\, .
$$
\end{proposition}

\begin{proof}
In view of the expression of $\Theta_{n-1}(x)$ in
(\ref{regular-thetas2}), we consider curves given by
$$
c=\frac{v}{1+tv},
$$
where $c>0$ is a constant, and $t\in [\ell_0,r_0]$. Note that
these curves are monotonically increasing on $[\ell_0,r_0]$, and
that the curve given by $v=\frac{c_1}{1-c_1t}$ lies ``above'' the
curve given by $v=\frac{c_2}{1-c_2t}$, if and only if $c_1>c_2$.

Now let ${\mathcal L}_{\alpha}$ be defined as the positive largest
$c$ for which the curve $c=\frac{v}{1+tv}$ lies in
$\Omega_{\alpha}$ for $t\in [\ell_0,r_0]$, i.e.,
$$
{\mathcal L}_{\alpha}=\max \left\{ c>0:\, \left(
t,\frac{c}{1-ct}\right)\in \Omega_{\alpha}, \text{ for all } t\in
[\ell_0,r_0]\right\} .
$$
It follows from Theorem~\ref{distribution-orbit} that for all
$z\leq {\mathcal L}_{\alpha}$, and for almost all
$G_q$-irrationals $x$ one has that
$$
\lim_{N\to \infty} \frac{1}{N} \# \left\{ n:1\leq n \leq N,\,
\Theta_n(x) < z \right\}  = \int_{\ell_0}^{r_0}\left(
\int_0^{\frac{z}{1-zt}} g_{\alpha}(t,v)\, {\rm d}v\right) {\rm
d}t= \lambda C_{q,\alpha}z\, ,
$$
where $C_{q,\alpha}$ is the normalizing constant of the invariant
measure (which has density $g_{\alpha}$). So we are left to show
that ${\mathcal L}_{\alpha} = \min\left\{
\frac{\lambda}{\lambda+2},\frac{\lambda
(2-\alpha\lambda^2)}{4-\lambda^2}\right\}$.

In the even case we can discern three cases: $\alpha=1/2$,
$1/2<\alpha<1/\lambda$, and $\alpha=1/\lambda$. Note that the
first case has been dealt with in \cite{[BKS]}; in case
$\alpha=1/2$ one has that ${\mathcal
L}_{\alpha}=\frac{\lambda}{\lambda +2}$.\smallskip\

In case $1/2<\alpha<1/\lambda$, first note that the curve
$c_1=\frac{v}{1+tv}$ goes through
$(r_1,H_1)=(\frac{1}{\alpha\lambda}-\lambda, \frac{1}{\lambda
+1})$ if and only if $c_1=\frac{\alpha\lambda}{\alpha\lambda +1}$.
Since in this case
$$
\frac{\alpha\lambda}{\alpha\lambda +1} < \frac{1}{2} <
\frac{1}{\lambda}=H_2,
$$
and the curve $c_1=\frac{v}{1+tv}$ is monotonically increasing on
$[\ell_0,r_0]$, we immediately find that this curve is in
$\Omega_{\alpha}$ for $t\in [\ell_0,0]$, yielding that $\mathcal
L_{\alpha}\leq \frac{\alpha\lambda}{\alpha\lambda +1}$.

From Theorem~\ref{thm:evenorder} we see that
$\ell_{p-2}<0<\ell_{p-1}$. Setting
$$
c_2=\frac{H_{2p-2}}{1+\ell_{p-1}H_{2p-2}}=\frac{\lambda
(2-\alpha\lambda^2)}{4-\lambda^2},\quad
c_3=\frac{H_{2p-1}}{1+r_0H_{p-1}}=\frac{1}{1+\alpha\lambda},
$$
it follows from Theorem~\ref{thm:even case nat ext} (see also
Figure~\ref{fig:nat ext q=6} for $q=6$) that ${\mathcal
L}_{\alpha}=\min \{ c_1,c_2\}$, since $c_1<c_3$. For $q$ (and
therefore $\lambda$) fixed, and for $\alpha\in [1/2,1/\lambda ]$,
one easily shows that $c_1=c_1(\alpha )$ is a monotonically
increasing function of $\alpha$, with $c_1(1/2)=
\frac{\lambda}{\lambda +2}$, and $c_1(1/\lambda )=1/2$, while
$c_2=c_2(\alpha )$ is a line with slope
$-\lambda^3/(4-\lambda^2)$. Since $c_2(1/2)=\lambda /2>1/2>\lambda
/(\lambda +2)$, and $c_2(1/\lambda )=\lambda /(\lambda
+2)=c_1(1/2)<c_1(1/\lambda )=1/2$, we find for $1/2<\alpha
<1/\lambda$ that
$$
{\mathcal L}_{\alpha} = \min\left\{
\frac{\lambda}{\lambda+2},\frac{\lambda
(2-\alpha\lambda^2)}{4-\lambda^2}\right\}.
$$
In case $\alpha =1/\lambda$, the point $(1,\lambda /2)$ yields
$c=\lambda /(\lambda +2)$. Since the curve $c=\frac{v}{1+tv}$ is
monotonically increasing on $[\ell_0,r_0]$, and from the fact that
for $t=0$ we have that $v=c=\lambda /(\lambda +2)<1/\lambda$,
where $1/\lambda$ is the ``smallest height'' of $\Omega_{\alpha}$,
we find that $$ {\mathcal L}_{1/\lambda}= \frac{\lambda}{\lambda
+2}.
$$
This proves the theorem.
\end{proof}

\begin{remark0} We only deal with the even case in
Proposition~\ref{lenstra}; a result for the odd case is obtained
similarly, but has a more involved
expression.\end{remark0}\smallskip\

\noindent \textbf{Acknowledgements.} We thank the referee for many
valuable remarks, which improved considerably the quality of the
presentation of this paper.


\end{document}